\documentclass[10pt]{article}
\usepackage[latin1]{inputenc}
\usepackage[english]{babel}
\usepackage[T1]{fontenc}
\usepackage{subfigure}
\usepackage{amssymb}
\usepackage{listings}   
\usepackage{amsmath}
\usepackage{pifont}
\usepackage{mathtools}
\usepackage{float}
\usepackage{listings}
\usepackage{graphicx}
\usepackage{hyperref}
\usepackage{caption}
\usepackage{multirow}
\usepackage[latin1]{inputenc}
\usepackage{framed} 
\usepackage{xcolor}
\usepackage{moreverb}
\usepackage{siunitx}
\usepackage{tabularx}
\usepackage{mathtools}
\colorlet{shadecolor}{gray!40}

\usepackage{authblk}
\usepackage{amsthm}
\newtheorem*{remark}{Remark}

\def\d{\,\,\textrm{d}}

\renewcommand{\epsilon}{\varepsilon}
\renewcommand{\phi}{\varphi}
\newcommand{\grad}{\nabla}
\renewcommand{\div}{\mbox{div\,}}

\usepackage{color, colortbl}
\definecolor{myred}{rgb}{0.93,0.76,0.75}

\DeclarePairedDelimiter{\abs}{\lvert}{\rvert}

\makeatletter
\let\@fnsymbol\@arabic
\makeatother

\renewcommand{\vec}[1]{\mbox{\boldmath $#1$}}

\newcommand*{\TitleFont}{%
      \usefont{\encodingdefault}{\rmdefault}{b}{n}%
      \fontsize{16}{20}%
      \selectfont}

\title{\TitleFont \textbf{A computational model of open-irrigated radiofrequency catheter ablation accounting for mechanical properties of the cardiac tissue}}
\author[1,*]{Argyrios Petras}

\author[1,2]{Massimiliano Leoni}

\author[3]{Jose M. Guerra}

\author[1,2]{Johan Jansson}

\author[1]{Luca Gerardo-Giorda}

\affil[1]{BCAM - Basque Center for Applied Mathematics, Bilbao, Spain} 
\affil[2]{Department of Computational Science and Technology, KTH Royal Institute of Technology, Stockholm, Sweden}
\affil[3]{Department of Cardiology, Hospital de la Santa Creu i Sant Pau, Barcelona, Spain}

\affil[*]{Alameda de Mazarredo 14, 48009 Bilbao, Bizkaia, Spain, apetras@bcamath.org}

\date{October 4th, 2018}

\begin{document}
\maketitle

\begin{abstract}
Radiofrequency catheter ablation (RFCA) is an effective treatment for cardiac arrhythmias. 
Although generally safe, it is not completely exempt from the risk of complications. 
The great flexibility of computational models can be a major asset in optimizing interventional strategies, if they can produce sufficiently precise estimations of the generated 
lesion for a given ablation protocol. This requires an accurate description of the catheter tip and the cardiac tissue. 
In particular, the deformation of the tissue under the catheter pressure during the ablation is an important aspect that is overlooked in the existing literature, that resorts to a 
sharp insertion of the catheter into an undeformed geometry. As the lesion size depends on the power dissipated in the tissue, and the latter depends on the percentage of the electrode surface 
in contact with the tissue itself, the sharp insertion geometry has the tendency to overestimate the lesion obtained, especially when a larger force is applied to the catheter. 
In this paper we introduce a full 3D computational model that takes into account the tissue elasticity, and is able to capture the tissue deformation 
and realistic power dissipation in the tissue. Numerical results in FEniCS-HPC are provided to validate the model against experimental data, and to compare the lesions obtained 
with the new model and with the classical ones featuring a sharp electrode insertion in the tissue.

\end{abstract}
\smallskip
\noindent \textbf{Keywords:} radiofrequency ablation; open-irrigated catheter; elastic tissue deformation; finite elements

\section{Introduction}
Radiofrequency catheter ablation (RFCA) is a minimally invasive therapy
widely used for the treatment of various types of cardiac arrhythmias
\cite{huang2014catheter,Guerra2013}. Typically, for endocardial
radiofrequency ablations (RFA), a catheter is advanced into the cardiac
chamber from the groin of the patient through a blood vessel.
Radiofrequency (RF) current at \SI{500}{\kHz} is delivered to the
arrhythmogenic tissue via the electrodes on the catheter,
producing resistive electrical heating in the neighborhood of the electrode. Conduction then propagates the heating through the rest of the tissue. 
 At a target temperature, usually of \SI{50}{\celsius}, irreversible damage is inflicted to the
tissue and a lesion is formed \cite{huang2014catheter}.

RFA is considered an effective and safe treatment for cardiac
arrhythmias, however a number of life threatening complications may
occur.  Excessive heating of the blood near the electrode reaching around 
\SI{80}{\celsius} leads to the denaturation of the blood proteins and the
formation of a coagulum at the tip of the catheter
\cite{huang2014catheter,Demolin2002}.  In addition, extreme heating of
the tissue at \SI{100}{\celsius} leads to steam formation and audible steam pops occur
\cite{haines1990observations, Koruth2013, Chik2015}, which could result in
myocardial tear or tamponade \cite{huang2014catheter}.

Several computational models have been developed to describe the
biophysics of the cardiac RFA treatment process. Some use an axisymmetric
approach and develop 2-D geometries, on which they employ Penne's bioheat equation
with convective boundary conditions for the cooling effect of the blood
flow, coupled with a quasi-static electrical potential equation
\cite{labonte1994numerical,shimko2000radio,tungjitkusolmun2000thermal,berjano2004thermal}.
To model the effect of the coolant liquid of open irrigated catheters,
these models typically fix the temperature of the tip electrode
\cite{irastorza2018thermal,gonzalez2018should}.

Three dimensional cardiac RFA models that use the same mathematical
equations as the 2-D axisymmetric models are also available
\cite{panescu1995three,berjano2006theoretical,Gallagher2011}. For open
irrigated catheters, some models include the interaction of the saline
inflow with a layer of saline in epicardial modelling
\cite{gopalakrishnan2002mathematical} or with the surrounding blood in
endocardial RFA \cite{gonzalez2016computational}. Typically, a sharp
insertion of the catheter in the tissue is considered, disregarding its
elastic deformation. However the need for the
incorporation of mechanical deformation in RFA modelling was hinted in a previous work
\cite{Gallagher2011}.  Indeed, the amount of power delivered to the
tissue depends on the surface of contact between the electrode and the
tissue itself \cite{wittkampf2006rf}. Sharp insertion produces larger surfaces of
contact, and is likely to overestimate the actual effect of the RFA
treatment. The first RFA model that includes elastic deformation of
the tissue uses profiles that are extracted from X-ray scanning before
the ablation process \cite{cao2002fem}. No other work was found in our
literature review that includes the tissue deformation due to the contact
with the catheter.

We propose a computational model that includes both physical and
mechanical properties of the tissue. We consider in this
work a 6-holes open electrode with a hemispherical tip, placed vertically on the cardiac tissue. 
We model the elastic deformation of the tissue in contact with the tip by
solving the axisymmetric Boussinesq problem for a spherical punch. We use a
modified version of Penne's bioheat equation to model the thermal effect
of RFA and a quasi-static electrical potential equation. Constant-power
ablation is obtained by augmenting the system with a power constraint
equation. The model also features the incompressible Navier-Stokes
equations to describe the flow of blood and
the saline flow from the open-irrigated electrode. We
validate our proposed computational model against a set of suitably
defined \emph{in-vitro} experiments performed at J.M. Guerra's lab at the
Hospital de la Santa Creu i San Pau in Barcelona.  We investigate the effect
of tissue elasticity in the computational RFA models by comparing the
computational lesion dimensions of the deformed tissue case against the
undeformed one that is typically used in state-of-the-art models.

The paper unfolds as follows. Section~\ref{sec_computationalModel} describes our mathematical
model, and in Section~\ref{sec_calibration} we describe its calibration. In Section~\ref{sec_results} we validate the proposed model against
a set of \emph{in-vitro} experiments: we describe the experimental setup
and discuss the results of the proposed computational model. We also
compare the results of our model with the ones obtained with a geometry
featuring a sharp insertion of the electrode. Finally, in Section~\ref{sec_discussion} we discuss our findings and
the limitations of the proposed model and suggest some ways to overcome them.

\section{Computational model}\label{sec_computationalModel}

Our mathematical model consists of a system of partial differential
equations that describe the evolution of the temperature, the flow of
blood and the cooling saline solution, as well as the electrical field
generated by the RFA procedure. In this section, we expand on
its main characteristics: the geometry, the governing equations, and the method to assess the computational lesion.

\subsection{Geometry}\label{sub_geometry}
The computational domain that we consider in this paper is based on an \emph{in-vitro} experimental setup, which is similar to the one
presented by Guerra et al. \cite{guerra2013effects}. We construct a box of dimensions
\SI{80 x 80 x 80}{\mm} that consists of the blood chamber
\SI{80 x 80 x 40}{\mm}, the cardiac tissue \SI{80 x 80 x 20}{\mm}, and a
board \SI{80 x 80 x 20}{\mm} underneath the tissue. The board separates the tissue from the dispersive electrode, that is located at the bottom of the computational domain.
This is an important aspect, matching the settings in the experiments used for validation, where the tissue is not in direct contact with the dispersive electrode, but is actually placed on a physical
board of methacrylate. Moreover, in clinical RFA, the dispersive electrode is placed on a patch on the back (or thigh) of the patient.
Having the dispersive electrode in direct contact with the tissue would lead to an excessive simplification of the model.
Also, the presence of the board allows us to match the impedance of the
RFA system (provided by the machine) without modifying the electrical conductivity of the tissue.

\begin{figure}
    \centering
    \includegraphics[width=0.6\textwidth]{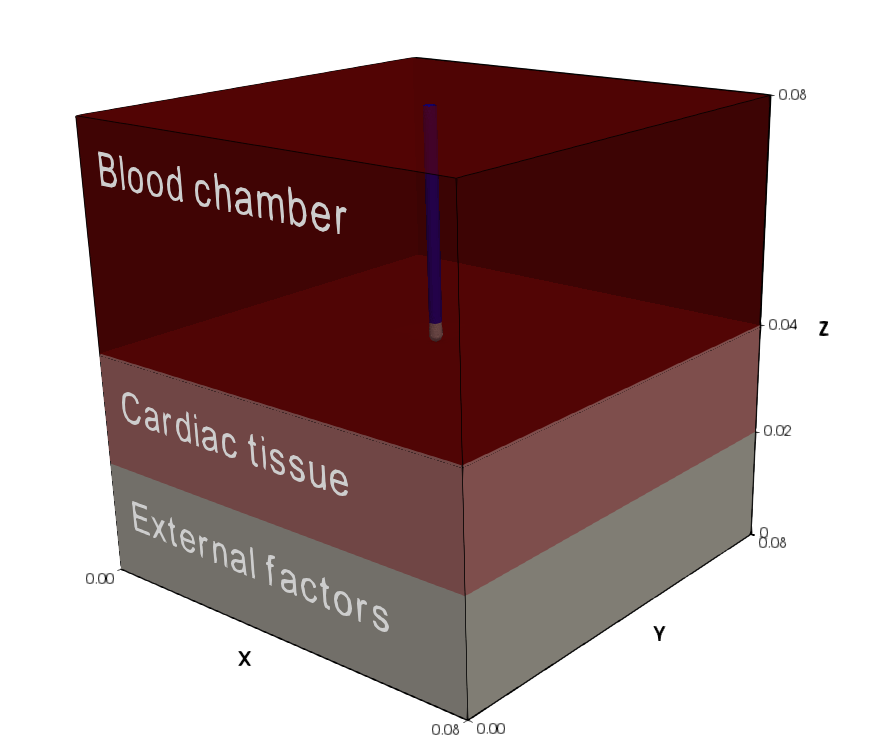} \ 
    \includegraphics[width=0.35\textwidth]{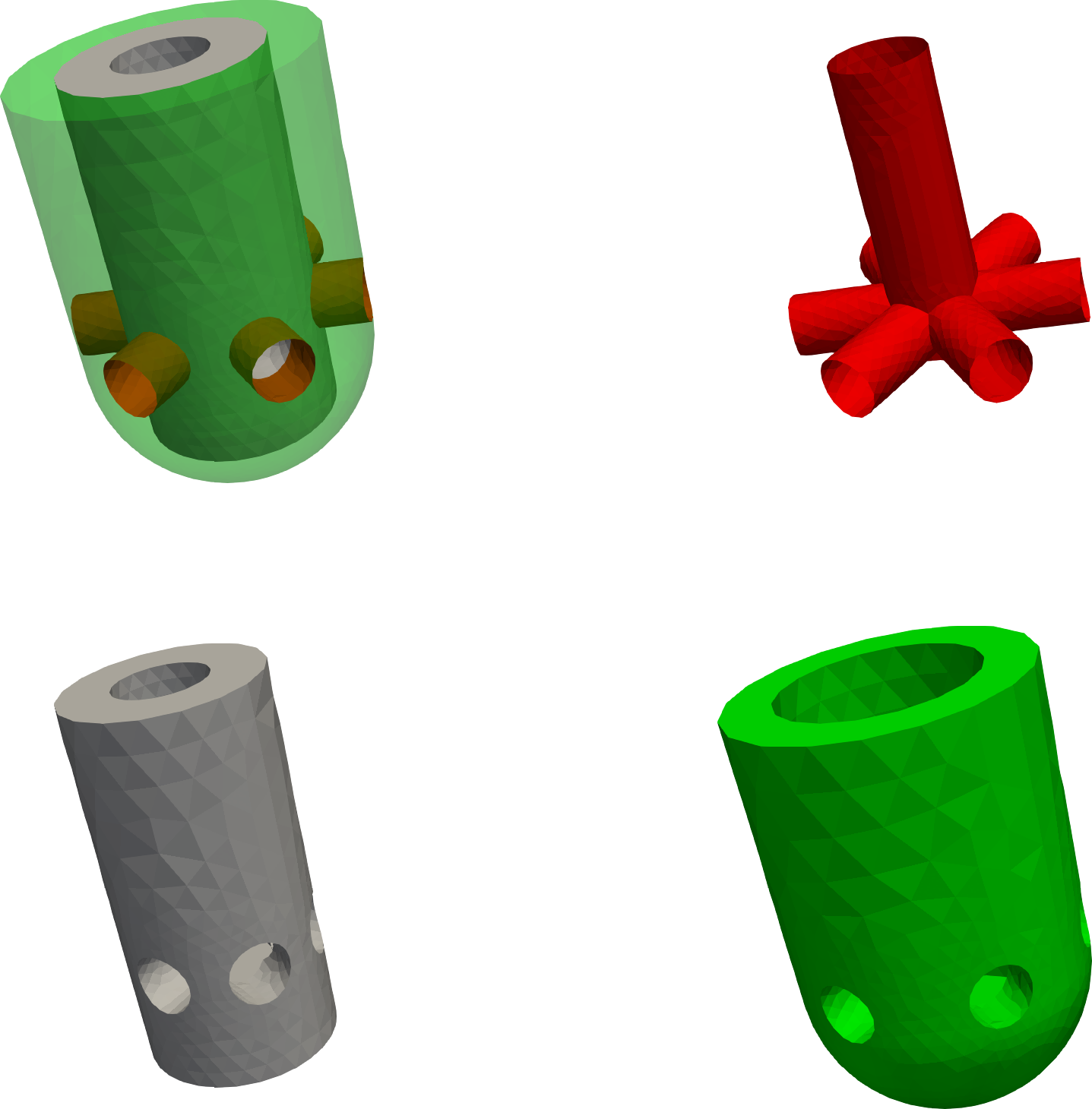}
   \caption{Left: The full computational geometry. Right: Detailed view of
  the catheter tip and its components. Top left to bottom right: the
  assembled catheter tip, the saline pipes, the thermistor, and the
  electrode.}
  \label{fig_compGeometry}
\end{figure}

We place a catheter perpendicular to the tissue at the center of
the box inside the blood chamber, causing a deformation that depends on
the catheter-tissue contact force (see Section~\ref{sub_elasticity}).
Figure~\ref{fig_compGeometry} (left) shows the described computational geometry.
We consider a 6-holes open-irrigated electrode of diameter \SI{2.33}{\mm}
and length \SI{3.5}{\mm} with a hemispherical tip. The electrode contains
a thermistor of diameter \SI{1.54}{\mm} and length \SI{3}{\mm}. The
holes, of diameter \SI{0.5}{\mm}, are connected to an inner channel of
diameter \SI{0.73}{\mm} inside the thermistor, which is connected to the
catheter body and allows the saline to flow in the blood chamber.
Figure~\ref{fig_compGeometry} (right) shows the geometry of the
open-irrigated electrode.

The computational domain consists of five different subdomains: the blood chamber ($\Omega_{blood}$), the cardiac tissue ($\Omega_{tissue}$), the external factors
board ($\Omega_{board}$), the electrode ($\Omega_{el}$) and the thermistor ($\Omega_{therm}$). In what follows, we will simply denote
$$
\Omega := \Omega_{blood} \cup \Omega_{tissue} \cup \Omega_{board} \cup \Omega_{el} \cup \Omega_{therm}.
$$
The catheter body and the inner tubes at the open-irrigated electrode are considered thermally and electrically insulated, and therefore insulation boundary
conditions are applied. Additional details on
the thermal and electrical boundary conditions are provided in
Sections~\ref{sub_potential} and \ref{sub_temperature}.

\subsection{Blood and saline flow}\label{sub_bloodFlow}
The open-irrigated electrode allows the cooling saline to flow out through the holes at its tip, and mix with the blood.
The comparable densities between saline and blood justifies the assumption of a perfect mix in the blood chamber. The flow dynamics is thus governed by the incompressible Navier-Stokes
equations:
\begin{equation}
  \begin{aligned}
    \frac{\partial \vec{u}}{\partial t} + \vec{u}\cdot
    \grad\vec{u} - \div \sigma(\vec{u},p) & = & \vec{0} \\
    \div\vec{u} &= & 0
  \end{aligned}
  \label{eq_navier_stokes}
\end{equation}
where $\vec{u}$ is the flow velocity, $p$ is the pressure scaled by the density $\rho$ and
\(\sigma(\vec{u},p) =  2 \mu\rho^{-1} \frac{\grad \vec{u} + \grad
\vec{u}^T}{2} - pI\) is the stress tensor, $\mu$ being the dynamic viscosity of the blood.

We impose zero pressure and stress as an outflow boundary condition on the right side of the domain $\{x=0.08\}$.
A constant inflow boundary condition is imposed on the left side of the domain $\{x=0\}$,
$$\vec{u} = \vec{u}_{in} = (u_b,0,0),$$
and on the electrode holes for the flow of the coolant liquid,
$$\vec{u} = u_s\vec{n}_e,$$
where $u_b$ and $u_s$ are the magnitudes of the velocity of the blood and the
coolant liquid respectively, while $\vec{n}_e$ is the unit normal vector
of the electrode pointing towards the blood chamber. All the remaining
boundaries, including the internal one separating the blood chamber from the
tissue, are equipped with no-slip boundary conditions. See Figure~\ref{fig_bcs} for a representation of all boundary conditions.

\subsection{Electrical field}\label{sub_potential}
The electrical field is generated by the electrode at the tip of the
catheter. At a frequency of \SI{500}{\kHz}, and over the distance
of interest, the biological medium can be considered totally
resistive and the electrical problem can be set in quasi-static form:
\begin{equation}
  \label{eq_pot}
  \div (\sigma(T)\grad\Phi)=0,
\end{equation}
where $\Phi$ is the electrical potential. We consider the electrical conductivity
$\sigma = \sigma(T)$ (not to be confused with the stress tensor
\(\sigma(\vec{u},p)\) introduced in \eqref{eq_navier_stokes}) to be
constant in
all regions but in the tissue, where we consider a linear
dependency on the temperature (see Section~\ref{sub_parameters}).

For a constant-power ablation mode, equation~(\ref{eq_pot}) is augmented with a constant-power constraint
\cite{Jiang2017}
\begin{equation}
  \label{eq_potentialMinimization}
  \int_\Omega \sigma(T)\grad\Phi\cdot\grad\Phi \d x = P,
\end{equation}
where $\Omega$ is the computational domain and $P$ is
the constant power imposed.

We set an initial potential of \SI{0}{\V} everywhere. We impose a zero potential
condition on the bottom of the domain to model the dispersive
electrode, while on the remaining external surfaces and the catheter
boundary we assign an insulation condition ($-\sigma \grad\Phi\cdot n = 0$).

\subsection{Tissue heating}\label{sub_temperature}
The application of an electrical potential at the tip electrode of the
catheter produces resistive heating at the cardiac tissue and the
surrounding blood. A modification of Penne's bioheat equation models both
the heating by the direct application of RF current and the conductive
heating \cite{gonzalez2016computational} and reads
\begin{equation}
  \frac{\partial (\rho h)}{\partial t} - \div(k(T)\grad T) = q - \rho
  c(T)\vec{u}\cdot\grad T + Q_m - Q_p,
\end{equation}
where $\rho$ is the density, $h$ is the specific enthalpy, $c(T)$ is the
specific heat, $k(T)$ is the thermal conductivity of the medium, $q =
\sigma(T)|\nabla \Phi|^2$ is the distributed heat source from the
electrical field, $Q_m$ is the metabolic heat generation and $Q_p$ is the
heat loss due to the blood perfusion. For endocardiac RFA, the last two
terms in this equation can be omitted for short ablation times
\cite{berjano2006theoretical,gonzalez2016computational,perez2017numerical}.

In biological tissues, the specific enthalpy $h$ and the temperature $T$
are related as follows:
\[
  \frac{\partial (\rho h)}{\partial t} = \psi(T)\frac{\partial T}{\partial t},
\]
where $\psi(T) = \rho c(T)$.
Finally, we assume the thermal conductivity $k(T)$ and the specific heat
$c(T)$ to be piecewise constant, with different values in different
regions of the domain, except in the tissue, where we consider a linear
dependence on the temperature (see Section~\ref{sub_parameters}).

The final form of the bioheat equation is then
\begin{equation}
  \rho c(T)\left(\frac{\partial T}{\partial t} + \vec{u}\cdot \grad
      T\right) - \div(k(T)\grad T) =
      \sigma(T)|\grad\Phi|^2.
  \label{eq_bioheatFinal}
\end{equation}

\subsection{Elasticity}\label{sub_elasticity}
During the ablation procedure, the electrode at the tip of the catheter applies a
contact force to the cardiac tissue, resulting in its deformation.
For a catheter with a hemispherical tip such as the one we consider in this paper
(see Figure~\ref{fig_compGeometry}), the problem is equivalent to a
spherical indenter on an elastic half-space in contact mechanics. There are
many different approaches for this contact mechanics problem, including
Hertzian theory for shallow indentations
\cite{asaro2006mechanics,popov2010contact}, a modified Hertzian solution
\cite{yoffe1984modified} and the axisymmetric Boussinesq problem
\cite{sneddon1965relation}. We use the axisymmetric Boussinesq problem approach as it provides a general framework for different profiles, which will allow us to consider also different catheter
tips in the future. In this work we consider a spherical indenter profile.

\begin{figure}[h!]
  \centering
  \includegraphics[width=0.7\textwidth]{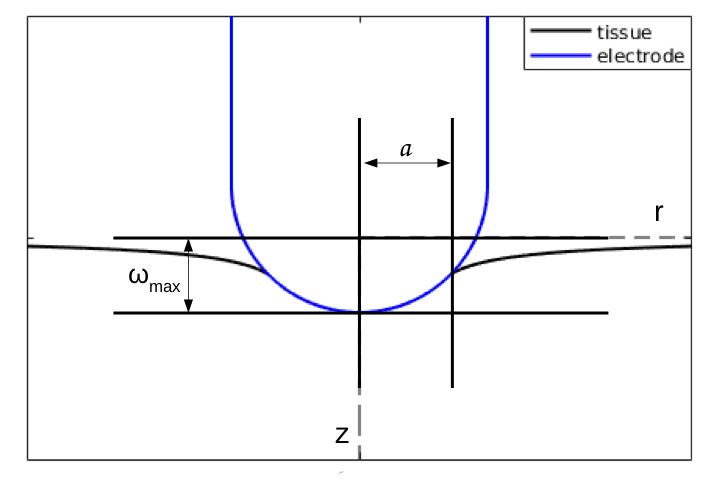}
  \caption{A cross-section of the electrode-tissue contact.}
  \label{fig_crossSection}
\end{figure}

Following \cite{sneddon1965relation}, due to the radial symmetry of the
contact problem, consider the two dimensional problem on the vertical
$z$-axis and the horizontal $r$-axis, which represents the distance from the center of
the sphere. The deformed elastic tissue follows the spherical indenter from
the center of the sphere, where the maximum vertical displacement
$\omega_{\max}$ occurs, up to the contact radius $a$. After the contact
point, the elastic tissue deformation is smoothly decreasing as the distance
$r$ increases. Figure~\ref{fig_crossSection} shows the deformation of
the tissue, as described.

Thus, for the case $r\leq a$, the elastic tissue follows the sphere, i.e.\
$$r^2 + [z+(R-\omega_{\max})]^2 = R^2,$$
where $R$ is the radius of the sphere. The positive solution of this system
for $z$ gives the displacement of the interface $\omega(r)$ as a function
of the distance $r$,
\begin{equation}
  \omega(r) = \omega_{\max} - f(r/a),
\end{equation}
where $f(t) = R - \sqrt{R^2-a^2t^2}$. The function $f$ is the profile of the
spherical indenter that is used in the rest of the computations.

The maximum depth is given by the formula
\begin{equation}
  \label{eq_depth}
  \omega_{\max} = \frac{a}{2}\log\left(\frac{R+a}{R-a}\right),
\end{equation}
and the total force is
\begin{equation}
  \label{eq_force}
  F = \frac{G}{1-\nu}\left((a^2+R^2)\log\left(\frac{R+a}{R-a}\right) - 2aR\right),
\end{equation}
where $\nu$ is the Poisson's ratio and
$$G = \frac{E}{2(1+\nu)}$$
is the shear modulus of the tissue, $E$ being the Young's modulus of elasticity.

For $r > a$, the deformation can be calculated by evaluating the integral
\begin{equation}
  \omega(r) = \int_0^1\frac{\chi(t)}{(r/a)^2-t^2}\;dt,
\end{equation}
where
$$\chi(t) = \frac{2\omega_{\max}}{\pi} - \frac{at}{\pi}\log\left(\frac{R+at}{R-at}\right).$$

Thus, given the contact force $F$, the contact radius $a$ can be
calculated by Equation~(\ref{eq_force}) using optimization
techniques. The remaining quantities can be calculated using the corresponding
formulas above. The displacement of the interface $\omega(r)$ is calculated using
numerical integration techniques for $r>a$.

\subsection{Model summary}
\label{sub_boundary_initial}
We consider a 3D computational geometry, as described in Section~\ref{sub_geometry}, and we compute the vertical displacement of the surface of the tissue for a given
catheter contact force (see Section~\ref{sub_elasticity}).

\begin{figure}[h!]
  \centering
  \includegraphics[width=0.9\textwidth]{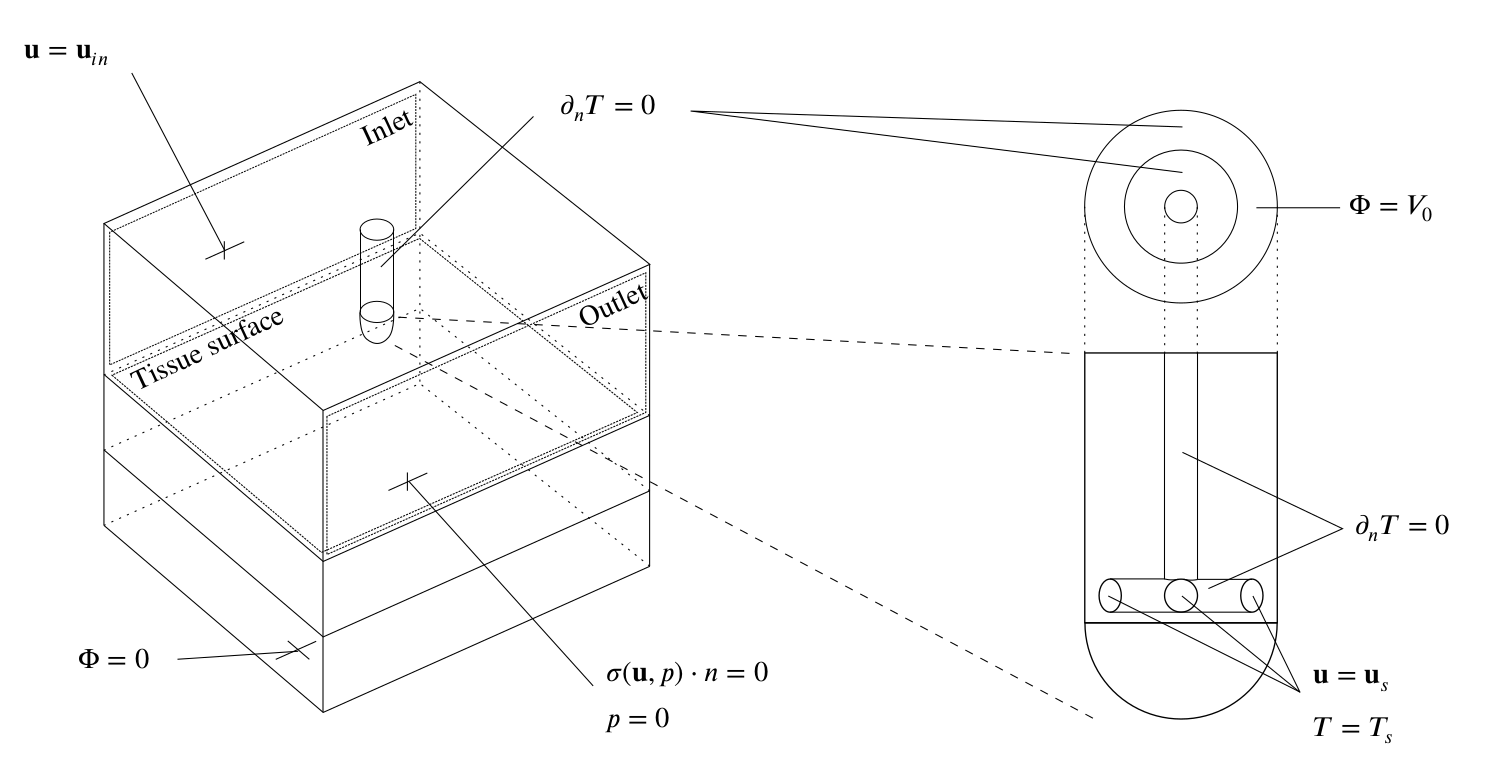}
  \caption{Boundary conditions that differ from zero
  velocity (\(\vec{u} = \vec{0}\)), zero current flux ($\sigma \grad\Phi\cdot n = 0$) and body temperature (\(T = T_b\)).}
  \label{fig_bcs}
\end{figure}

We summarize the mathematical model in the following system of coupled PDEs:
\begin{subequations}
  \begin{align}
    \frac{\partial \vec{u}}{\partial t} + \vec{u}\cdot \grad\vec{u}
      - \div\sigma(\vec{u},p) & = \vec{0}  & \mbox{in }\  \Omega_{blood} \times (0,T),\label{eq_NS}\\
    \div\vec{u} & = 0 & \mbox{in }\  \Omega_{blood} \times (0,T),\label{eq_NSinc}\\
    \rho c(T)\left(\frac{\partial T}{\partial t} + \vec{u}\cdot \grad
      T\right) - \div(k(T)\grad T) & =
      \sigma(T)|\grad\Phi|^2  &  \mbox{in }\  \Omega \times (0,T), \label{eq_bioheat}\\
    \div(\sigma(T)\grad\Phi) & = 0 & \mbox{in }\  \Omega \times (0,T).\label{eq_potential}
  \end{align}
  \label{eq_fullSystem}
\end{subequations}
The above equations are equipped with a number of different boundary
conditions, applied to the various parts of our computational geometry. Figure~\ref{fig_bcs} collects the boundary conditions of \eqref{eq_NS}-\eqref{eq_potential}: unless
stated otherwise, we impose zero velocity (\(\vec{u} = \vec{0}\)), zero current flux ($\sigma \grad\Phi\cdot n = 0$) and body temperature (\(T = T_b = \SI{37}{\celsius}\)) on
all the surfaces of the boundary.
We impose zero velocity on the tissue surface, which is part of the
boundary of $\Omega_{blood}$. The values of the saline inlet $\vec{u}_s$
are given by the RF protocol under consideration.  The homogeneous
Dirichlet boundary condition for the potential $\Phi = 0$ on the bottom
surface models the dispersive electrode, while the Dirichlet condition
$\Phi = V_0$ on the upper boundary of the electrode is discussed in the
upcoming Sections \ref{sub_board_calibration} and
\ref{sub_dirichletBC_update}. We set the initial conditions to \(p = 0\),
\(\vec{u} = \vec{0}\), \(\Phi = 0\) and \(T = T_b\) in the whole domain.

\subsection{Lesion assessment}
The prediction of the lesion size is important to ensure the
effectiveness of RF treatment.  During the ablation process, RF current
flows through the myocardium and, via Joule heating, a rise in the tissue
temperature is produced in a small area around the electrode tip.  The
tissue is irreversibly damaged and a thermal lesion is formed when the
temperature rises above \SI{50}{\celsius} \cite{huang2014catheter}. In
light of the above consideration, we assume that a lesion is fully formed
at the isotherm contour of \SI{50}{\celsius}.  In order to assess the
lesion characteristics, typical quantities of interest are its depth (D),
width (W), depth at which the maximal width occurs (DW), as well as its
volume (V). Also, the area of the lesion on the tissue surface (S) is
important to study the possible formation of thrombi. The characteristic
dimensions we consider in this paper to validate the computational model
are described in Section \ref{sub_expSetup}, and collected in Figure~\ref{fig_lesionDimensions} therein.

\section{Calibration of the model}
\label{sec_calibration}

\subsection{Model parameters}\label{sub_parameters}
The parameters considered in this paper are drawn from the literature and
are summarized in Table~\ref{tab_parameters}. Specifically, the
physical properties of the blood are obtained from the virtual population
database from the IT IS foundation \cite{hasgall2015database}. The
electrode and thermistor physical properties were found in the literature
\cite{Gallagher2011,gonzalez2016computational}. We assume the electrical
conductivity of the external factors board $\sigma_b$ to be constant and
we tune it to match the initial impedance of the experimental setup, as detailed in Section \ref{sub_board_calibration}.
We choose the remaining parameters for the board to be the same as the
initial tissue state; they appear to have no impact on the results of our numerical
simulations.

\begin{table}[h]
  \centering
  \begin{tabular}{lccccc}
    \hline
    Parameter & Blood & Tissue & Electrode & Thermistor & Board\\ \hline
    $\rho$ {\scriptsize(\si{\kg\per\m\cubed})} & 1050 & 1076 & 21500 & 32 & 1076\\
    $c$ {\scriptsize(\si{\J\per\kg\per\K})} & 3617 & $c_0 = $ 3017 & 132 & 835 & 3017\\
    $k$ {\scriptsize(\si{\W\per\m\per\K})} & 0.52 & $k_0 = $ 0.518 & 71 & 0.038 & 0.518\\
    $\sigma$ {\scriptsize(\si{\siemens\per\m})} & 0.748 & $\sigma_0 = $ 0.54 & $4.6\times10^6$ & $10^{-5}$ & $\sigma_{b}$\\
    $\mu\rho^{-1}$ {\scriptsize(\si{\m\squared\per\s})} & $2.52\times10^{-6}$ & - & - & - & -\\
    $\nu$ {\scriptsize(-)} & - & 0.499 & - & - & -\\
    $E$ {\scriptsize(\si{\kg\per\m\per\s\squared})} & - & $75\times10^3$ & - & - & -\\ \hline
  \end{tabular}
  \caption{The parameters we considered in our proposed computational model.}
  \label{tab_parameters}
\end{table}

Considering that in our \emph{in-vitro} experiments RFA was performed
on a porcine myocardium, we used the corresponding physical properties
\cite{tsai2002vivo,duck2013physical,bhavaraju1999thermophysical}.
We use a linear temperature-dependence of specific heat, thermal and
electrical conductivities of the tissue according to the aforementioned literature:

$$
c(T)  = c_0(1-0.0042(T-37)),\qquad
k(T)  =  k_0(1-0.0005(T-37)),\qquad
\sigma(T) = \sigma_0(1+0.015(T-37))
$$
where $c_0$, $k_0$ and $\sigma_0$ are the specific heat, the thermal
and the electrical conductivities at body temperature
(see Table~\ref{tab_parameters}). Finally, since the experiments were
performed in an \emph{in-vitro} condition, it is reasonable to assume that the cardiac tissue got stiffer during preparation. 
We thus used in our simulations the Young's modulus of the cardiac tissue at systole \cite{urban2013measurement}, which is within the range of measurements in an \emph{in-vitro} experimental setup found in \cite{ramadan2017standardized}.

\subsubsection{Tissue power delivery calibration}\label{sub_pow_calibration}

In our model, we set the electrical conductivity of the board $\sigma_b$ to match the power dissipated in the tissue $P_{tissue}$, that we identify following the analysis of Wittkampf and
Nakagawa in \cite{wittkampf2006rf}: the percentage of the total power $P$  delivered to the tissue is given by the formula:
\begin{equation}\label{eq_tissuePowerDissipation}
P_{tissue} =  \frac{A_{tissue}\sigma_{tissue}}{A_{blood}\sigma_{blood} + A_{tissue}\sigma_{tissue}}\, P =: \alpha P,
\end{equation}
where $P$ is the total power set
by the machine in the experiment and $A_{blood}$ and $A_{tissue}$ is the contact area of
the electrode with the blood and the tissue respectively.

In our computational geometry (see Section~\ref{sub_geometry}) we chose a
hemispherical catheter tip. Thus, the surface area of the
catheter in contact with the tissue $A_{tissue}$ can be calculated as
follows
\[
  A_{tissue} =
  \left\{
    \begin{array}{ll}
      2\pi R h, & \text{if } h\leq R,\\
      2\pi R^2 + 2\pi R (h-R) - 6 R_h^2\cos^{-1}\left(1 + \frac{R - h}{R_h}\right), & \text{if } R< h \leq R+R_h\\
      2\pi R^2 + 2\pi R (h-R) - 6R_h^2\left(\frac{\pi}{2} + \cos^{-1}\left(2 +\frac{R -h}{R_h}\right)\right), & \text{if } R+R_h < h \leq R+2R_h\\
      2\pi R^2 + 2 \pi R (h - R) - 6\pi R_h^2, & \text{if } h\geq R+2R_h,
    \end{array}
  \right.
\]
where $h = \omega_{\max} - \omega(a)$ is the contact depth of the
catheter with the tissue (see Section~\ref{sub_elasticity}) and $R_h$ is
the radius of the irrigation holes on the catheter (see
Section~\ref{sub_geometry}). The blood surface area is calculated using
the formula
\[
A_{blood} = 2\pi R^2 + 2\pi R (h_e - R) - 6\pi R_h^2 - A_{tissue},
\]
where $h_e$ is the length of the electrode.

\subsubsection{External factors calibration}
\label{sub_board_calibration}

We included the external factors board in the model to represent the part
of the experimental setup, between the tissue and the dispersive electrode, that is not modelled.
We tune the external-factors board's conductivity $\sigma_b$ to match
the actual power delivered to the tissue as described in
Section~\ref{sec_computationalModel}.
This approach allows us to use as much of the information available from
the RFA system as possible, and provides us with great flexibility should
different types of tissue be considered, in a major
improvement with respect to standard approaches on the topic (including earlier works by some of the authors \cite{gonzalez2016computational}) that tuned the tissue conductivity to match the
overall system's initial impedance.

Given the total power of the ablation $P$ and the initial impedance of the system $R$, (both values available from the RFA system), the initial potential drop in the system
is given by Ohm's law as $V_0 = \sqrt{PR}$.
We can assume that the electrode cable is not dissipating energy, thus we can consider the potential drop to occur entirely between the cable-electrode junction and the dispersive electrode.
We thus impose a Dirichlet boundary condition on the upper boundary of
the electrode (see Figure~\ref{fig_bcs}), and we solve the potential
equation in the whole domain $\Omega$. The advantage of such a framework is that the interfaces between different materials are actually internal discontinuity surfaces of the domain: 
the equation itself takes care of the conductances drop across internal layers as a consequence of the solution's continuities \cite{nedelec2001acoustic}. In addition, modifications of the 
electrode material would be immediate.

The calibration of the external factors board conductivity is an iterative pre-processing step.
Given an initial guess of \(\sigma_b\), we solve equation~\eqref{eq_potential} with  $V_0=\sqrt{PR}$ as boundary condition, and suitably update the value of
$\sigma_b$, until the constraint
\begin{equation}
  \int_{\Omega_{tissue}} \sigma_0\abs{\grad\Phi}^2 \d x = P_{tissue}
  \label{eq_tissuePowerConstraint}
\end{equation}
is satisfied, $P_{tissue}$ being the one in
\eqref{eq_tissuePowerDissipation}. In the code developed to obtain the
results of Section~\ref{sec_results}, we implement a root-finding
algorithm via bisection search, that stops once the difference between
two subsequent iterations is smaller than \SI{0.01}{\W}. The calculated \(\sigma_b\) is kept constant throughout the simulation.

\subsubsection{Voltage calibration during the simulation}
\label{sub_dirichletBC_update}

In Section~\ref{sub_potential} no boundary conditions were specified for the potential equation at the upper part of the electrode, but we required the
solution to satisfy a constraint on the total dissipated power \eqref{eq_potentialMinimization}.
However, due to the difference in size between our computational domain and the actual ablation site, the total power dissipated in the computational setting
is expected to be smaller than the one provided by the machine. We assume the power dissipated outside our computational domain as well as the one in \(\Omega\) to be constant in time.
During the simulation we update the potential \(V_0\) on the upper interface of the electrode, so that the solution satisfies
\begin{equation}
 \label{eq_practical_constraint}
 \int_\Omega \sigma\abs{\grad\Phi}^2 \d x = P_0.
\end{equation}
In the above formula, $P_0$ is the total power dissipated in our computational domain at the beginning of the simulation with $\sigma_b$ computed as described in the previous section.
We assume a constant voltage on the surface where we impose the boundary condition: the electrode is made of conductive material and we know from standard theory of
electromagnetism that electrons can move freely on the surface of a conductor, tending to a configuration of equilibrium.

\section{Numerical Results}
\label{sec_results}
For the simulations in this work, we built a 3D computational geometry for each given catheter contact force, where the vertical displacement of the surface of the tissue
is computed by means of a Clenshaw-Curtis numerical integration method \cite{jones2016others}.

The resulting domain $\Omega$ is discretized by a tetrahedral unstructured mesh, with adaptive refinement in the neighborhood of the electrode.
The mesh is built with Salome \cite{ribes2007salome}. For all the
simulations presented in this work, the meshes feature approximately
\num{5000000} elements
with minimum and maximum element sizes of the order of magnitude of
\SI{5e-5}{\m} and \SI{1e-3}{\m}.
Following the instructions of the catheter producer, the saline inflow is set at $u_s =\SI[per-mode=symbol]{0.24}{\m\per\s}$, corresponding to an irrigation rate of $\SI[per-mode=symbol]{17}{\mL\per\min}$ subdivided among the 6 irrigation holes. In the model validation section, the blood inflow velocity is set to $u_b = \SI[per-mode=symbol]{0.5}{\m\per\s}$
to match our experimental setup described in the upcoming Section~\ref{sub_expSetup}.
In the comparison between elastic deformation and sharp insertion, two blood flow protocols are considered, as in the experiments by J. Guerra and his collaborators \cite{guerra2013effects}:
a High-Flow protocol (HF) with $u_b = \SI[per-mode=symbol]{0.5}{\m\per\s}$, and a Low-Flow protocol (LF), with  $u_b = \SI[per-mode=symbol]{0.1}{\m\per\s}$.

\subsection{Simulations protocol}\label{sub_simulation_protocol}

Computational modelling of ablation as well as other medical procedures and diagnostics in cardiovascular science represent significant scientific challenges, such as 
robust turbulent fluid-structure interaction (FSI) with contact and heating, control of the computational error in multiphysics problems, and complex software implementation.
For computational technology to be effective in overcoming these challenges, it
has to be {\bfseries efficient}: enabling solution of advanced problems with
available resources, {\bfseries reliable}: providing quantitative
error control or the computational solution in
a chosen output, as well as a guarantee of correctness of software
implementation, and {\bfseries accessible}: enabling scientists to
access, use and extend the simulation technology as a tool in their daily work.

We address these problems in the setting of a general stabilized
adaptive finite element methodology we denote Direct FEM Simulation
(DFS) \cite{hoffman2017computability}
realized in the automated massively parallel FEniCS-HPC open source
software framework
\cite{fenicsbook,fenicsppam}, taking the weak form of the partial
differential equation (PDE) as
input and automatically generating low-level source code for
assembling tensors and duality-based a posteriori error estimates and
indicators for adaptive error control with good scaling on supercomputers.
In this framework, we developed an RFA module in FEniCS-HPC explicitly for the present work to discretize and solve system \eqref{eq_NS}-\eqref{eq_potential}.

\paragraph{Navier-Stokes equations}
The Navier-Stokes equations are solved via the DFS method presented in \cite{hoffman2017computability, hoffman2013unicorn, jansson2018time}, consisting 
of a standard P1-P1 finite elements method with Galerkin Least-Squares stabilization and optional goal-oriented adaptive mesh refinement.
The Navier-Stokes system is discretized in time by a Crank-Nicholson scheme, and DFS solves the full nonlinear problem via successive
minimizations based on a fixed-point approach \cite{hoffman2017computability, jansson2018time}.

\paragraph{Bioheat equation}
The Bioheat equation, on the other hand, is solved with P1 elements and is advanced in time by a backward Euler scheme, where the coefficients are linearized around the previous time step.
Since for the set of parameters introduced in
Section~\ref{sub_parameters} the Bioheat equation is advection-dominated
in $\Omega_{blood}$, an SUPG (Streamline-Upwind Petrov-Galerkin)
stabilization term is added to the discrete equation\cite{quarteroni2008numerical}.

\paragraph{Potential equation}
The potential equation is solved with P1 elements.
The equation depends implicitly on time through the conductivity. As the simulation proceeds, the rise in the tissue temperature modifies its electrical conductivity, and the constraint on the total
dissipated power will eventually fail to be satisfied. We thus enforce constraint \eqref{eq_practical_constraint} at each time step when we solve the potential equation as follows.

Let $\Phi_1$ be the solution of
\begin{equation}
  \begin{aligned}
    \div(\sigma^{new} \grad\Phi_1) &= 0 && \text{ in } \Omega \\
    \Phi_1 &= V_0^{old} && \text{ on } \partial\Omega_{el}^{up} \\
    \Phi_1 &= 0 && \text{ on } \partial\Omega^{bottom} \\
    \sigma^{new} \grad\Phi_1 \cdot n&= 0 && \text{ on } \partial\Omega\backslash\{\partial\Omega_{el}^{up} \cup\partial\Omega^{bottom}\}
  \end{aligned}
  \label{eq_potentialFull}
\end{equation}
where $V_0^{old}$ and $ \sigma^{new}$ are the previous value of the potential on the upper part of the electrode and the new value of the electrical conductivity, respectively.
The power dissipated by $\Phi_1$ is
\[
  P_1 = \int_\Omega \sigma^{new} \abs{\grad\Phi_1}^2 \d x.
\]
If $|P_1-P_0| < $ \SI{0.01}{\W}, we choose $\Phi^{new}=\Phi_1$ as new solution.

On the other hand, if $P_1$ does not fall within the tolerance, we proceed as follows. By linearity, the function \(\Phi_2 = \lambda\Phi_1\) is the solution of problem \eqref{eq_potentialFull} with Dirichlet boundary condition \(\Phi_2 = \lambda V_0^{old}\) on $\partial\Omega_{el}^{up}$ (Dirichlet boundary conditions scale linearly and the Neumann one is itself homogeneous).
The power dissipated in the domain by $\Phi_2$ is
\[
  P_2 =\int_\Omega \sigma^{new}  \abs{\grad\Phi_2}^2 \d x =
  \int_\Omega \sigma^{new}  \abs{\grad(\lambda\Phi_1)}^2 \d x = \lambda^2
  \int_\Omega \sigma^{new}  \abs{\grad\Phi_1}^2 \d x = \lambda^2 P_1.
\]
By requesting $\Phi_2$ to satisfy constraint \eqref{eq_practical_constraint} (namely, imposing $P_2=P_0$), we can solve for $\lambda=\sqrt{{P_0}/{P_1}}$, and the new potential will be
$$
\Phi^{new} = \sqrt{\frac{P_0}{P_1}}\ \Phi_1, \qquad V_0^{new} = \sqrt{\frac{P_0}{P_1}}\ V_0^{old}.
$$

\paragraph{Time stepping}
A time step of $\Delta t = \SI{0.01}{\s}$ has been used for the bioheat equation and a smaller, adaptively computed time step is used for the
Navier-Stokes equations. As a consequence of the different time steps employed, the two models are synchronized at the larger time step, where also the potential equation is solved.
Let $t^n= n \Delta t$ be the $n$-th time step of the bioheat equation, and let \((\vec{u}_h^n, p_h^n, \Phi_h^n, T_h^n)\) be the numerical solution at time \(t^n\). We can briefly summarize the
procedure to advance from $t^n$ to $t^{n+1}$ as follows.

\begin{enumerate}
\item Solve the Navier-Stokes equations in $(t^n,t^{n+1})$ with the adaptive time step and compute, in particular,
\(\vec{u}_h^{n+1}\).
\item Solve the potential equation for \(\Phi_h^{n+1}\) with coefficient $\sigma(T_h^n)$
\item Solve the bioheat equation for \(T_h^{n+1}\) using \(\vec{u}_h^{n+1}\) in the transport term, \(\Phi_h^{n+1}\) in the source term,  $\sigma(T_h^n)$, $c(T_h^n)$ and $k(T_h^n)$ as coefficients.
\end{enumerate}

We refer the reader to the supplementary material section for a more detailed description of the numerical discretization.
Finally, the lesion dimensions are calculated in a post-processing step in Paraview \cite{ayachit2015paraview} with self-developed Python-coded filters.

\medskip

\begin{remark}
Since the stress tensor \(\sigma(\vec{u},p)\) is independent of the
temperature, equations~\eqref{eq_NS} and \eqref{eq_NSinc} are in a
one-way coupling with the rest of the system.
Thus, already in our current situation, but even should the geometry change in
time due to a variable pressure during
the ablation (as it is the case in clinical RFA procedures), it would
still be possible to solve the Navier-Stokes equation independently as a
first step and then use the solution to solve
equations~\eqref{eq_potential} and \eqref{eq_bioheat}, potentially for a
variety of different protocols. However, we
chose to solve the Navier-Stokes equations alongside the potential and
bioheat ones during each simulation in order to maintain a design 
flexible enough to handle also more complex (i.e.~fully
coupled) situations, such as a stress tensor of the tissue that varies
with temperature.
\end{remark}

\subsection{\emph{In-vitro} experimental setup}\label{sub_expSetup}
Two \emph{in-vitro} sets of experiments using porcine myocardial tissue were
performed in an experimental setup similar to the one presented in \cite{Guerra2013}. The
porcine myocardial tissue was placed on a polymethyl methacrylate board
in a bloodbath. An iron arm attached to the board kept the
catheter perpendicular to the tissue and maintained the
contact force around the target  \SI{10}{\g}, \SI{20}{\g} and
\SI{40}{\g}, as measured by the catheter system. A pump attached to the iron arm generated a constant
blood flow of \SI[per-mode=symbol]{0.5}{\m\per\s}. The indifferent
electrode was placed at the bottom of the blood bath. A heater preserved
the blood temperature to \SI{37}{\celsius}.

The 6-holes open-irrigated TactiCath Quartz catheter by St.\ Jude Medical
Inc.\ was used in the experiment. A constant power delivery ablation
protocol of \SI{20}{\W} or \SI{35}{\W} was followed for all the
ablations and a saline flow rate of
\SI[per-mode=symbol]{17}{\mL\per\minute} was set.

\begin{figure}[h!]
  \centering
  \includegraphics[width=0.4\textwidth]{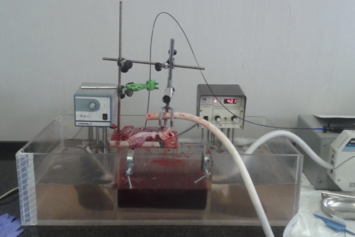} \qquad
  \includegraphics[width=0.4\textwidth]{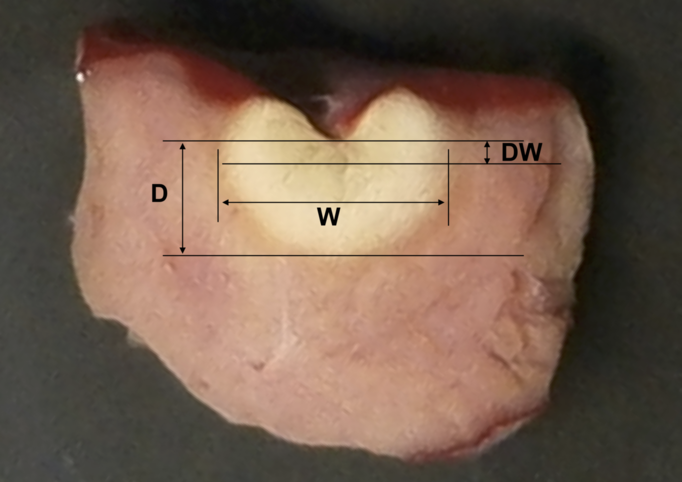}
  \caption{Left: the experimental setup. Right: lesion measurements: depth (D), width (W) and depth of maximum
  width (DW).}
  \label{fig_lesionDimensions}
\end{figure}

\begin{table}[h!]
  \centering
  \begin{tabular}{l|cccc}
    &\multicolumn{4}{c}{Target experimental setup ($F$, $P$, $T$)}\\ \hline
    & $(10^*\SI{}{\g},20\SI{}{\W},30\SI{}{\s})$ &
    $(10^*\SI{}{\g},35\SI{}{\W},30\SI{}{\s})$ & $(20^*\SI{}{\g},20\SI{}{\W},30\SI{}{\s})$ &
    $(40^*\SI{}{\g},20\SI{}{\W},30\SI{}{\s})$\\\hline
    $n$ & 8 & 5 & 9 & 2\\\hline
    $\bar{F}_{exp}$ {\scriptsize(\SI{}{\g})} & $11.2\;[11.1:11.5]$ & $11.9\;[11.3:12.0]$ & $19.9\;[19.6:22.2]$ & $41.0\;[40.5:41.5]$\\\hline
    D {\scriptsize(\SI{}{\mm})} & $2.3\;[1.8:2.6]$ & $4.4\;[4.3:4.6]$ & $3.1\;[2.8:3.5]$ & $3.3\;[3.0:3.5]$\\
    \hfill min & 1.4 & 4.2 & 2.2 & 3.0\\
    \hfill max & 3.1 & 4.6 & 4.6 & 3.5 \\\hline
    W {\scriptsize(\SI{}{\mm})} & $6.5\;[6.3:7.4]$ & $12.7\;[9.9:13.2]$ & $8.9\;[8.0:10.5]$ & $7.2\;[7.1:7.2]$\\
    \hfill min & 5.5 & 8.4 & 7.0 & 7.1\\
    \hfill max & 7.9 & 14.1 & 11.7 & 7.2\\\hline
    DW {\scriptsize(\SI{}{\mm})} & $0.4\;[0.4:0.7]$ & $0.4\;[-0.1:0.9]$ & $0.5\;[0.3:0.8]$ & $0.7\;[0.3:1.0]$\\
    \hfill min & 0.3 & $-$1.0 & $-$0.6 & 0.3\\
    \hfill max & 0.8 & 1.0 & 1.4 & 1.0\\\hline
    $R_0$ {\scriptsize(\SI{}{\ohm})} & $121\;[116:134]$ & $104\;[102:138]$ & $122\;[114:124]$ & $127\;[117:136]$\\
  \end{tabular}
  \caption{Summary of lesion dimensions for every experimental
  setup. The starred values $^*$ denote the mean nominal forces and
  include a standard deviation of \SI{1}{\g}. The values in the table are
  either absolute numbers or expressed in the format
  \emph{$Q_2$~[$Q_1$~:~$Q_3$]}, where $Q_2$ is the median and $Q_1$ and
  $Q_3$ are the first and the third quartiles.}
  \label{tab_lesionSummary}
\end{table}

After the radiofrequency delivery, the lesion dimensions were assessed by
examining their cross-section. The measured variables were the depth
(D), the width (W) and the depth of the maximum width (DW) as shown in
Figure~\ref{fig_lesionDimensions}. In standard computational models that feature a sharp insertion of the catheter in the tissue,
the depth and depth of maximum width are measured from the tissue surface. However, such measurement techniques have their limitations.
First, existing computational models do not take into account the irregular surface of the cardiac wall, and neither do we at this stage.
In addition, since lesions are assessed after physical manipulation, it is difficult to identify the actual position of the pre-ablation tissue surface.
By considering these limitations, we opted to use measurements that allow us to validate our numerical results through comparable variables (see Figure~\ref{fig_lesionDimensions} and
Figure~\ref{fig_numerical_lesion}, right).

A total of \num{31} ablations were performed for different target experimental
setups, using three control parameters: the contact force $F$, the
applied power $P$ and the time of ablation $T$. During the
experiments we observed that the target contact force was not
always achieved due to maneuvrability limitations and the effect of the blood inflow to the
system. Thus, we allow a standard deviation of \SI{1}{\g} as error of
the mean experimental force $F_{exp}$ around the
nominal force $F$ and we remove the observations that deviate three standard
deviations. Table~\ref{tab_lesionSummary} summarizes the
lesion dimensions and the initial impedance of the system $R_0$. The mean experimental contact force $\bar{F}_{exp}$ is
reported for the non-extreme cases. No complications (steam pops, char
formation etc.) were observed in our experiments.

\subsection{Model validation}
\label{sub_Validation}

\begin{figure}[h!]
  \centering
  \includegraphics[width=0.4\textwidth]{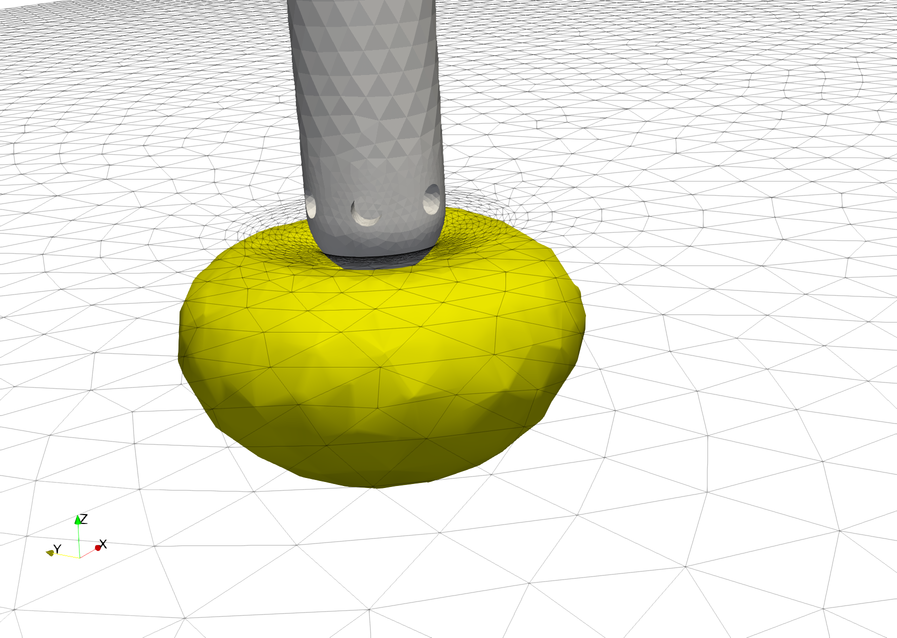} \
  \includegraphics[width=0.4\textwidth]{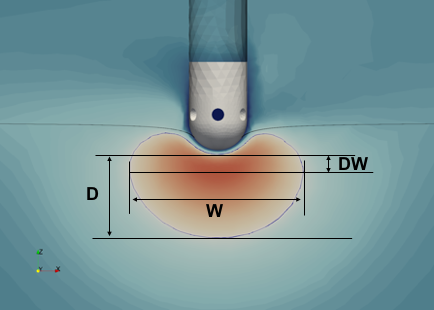}
  \caption{Left: the 3D computed lesion in the (20\SI{}{\g}, 20\SI{}{\W}, 30\SI{}{\s}) settings. Right: lesion measurements used for validation.}
  \label{fig_numerical_lesion}
\end{figure}

We validate our model by comparing the computational lesion dimensions against the \emph{in-vitro} experimental results in Section~\ref{sub_expSetup}.
We run simulations with parameters corresponding to the first 3 settings of the experiments. We discarded the (40g, 20W, 30s) configuration because only two
experiments were available, and one of the two images obtained from the
experiments was too blurry to get solid measurements.
As an example, in Figure~\ref{fig_numerical_lesion} we show the 3D lesion computed in the (20\SI{}{\g}, 20\SI{}{\W}, 30\SI{}{\s}) settings, and the corresponding lesion dimensions used for
validation. Given the uncertainty on the cross-section orientation (with respect to the blood flow direction) of the experimental lesions, these dimensions are assessed on the plane
 featuring the maximum width of the lesion.

In Table~\ref{tab_minMaxComparison} we collect the resulting lesion dimensions of our numerical simulations. We
observe that the depth $D$, the width $W$ and the depth of the maximum
width $DW$ of the computational lesion are within the range of the
corresponding experimental values for the case of $F=\SI{20}{\g}$. In the \SI{10}{\g} cases, the lesion width is underestimated by around $10\%$ for both cases. Furthermore, the depth of the maximum width is overestimated by nearly $50\%$ for the case of $P=\SI{35}{\W}$. 

\begin{table}[h!]
  \centering
  \begin{tabular}{l|ccc}
    &\multicolumn{3}{c}{Computational setup ($F$, $P$, $T$)}\\ \hline
    & $(10\SI{}{\g},20\SI{}{\W},30\SI{}{\s})$ &
    $(10\SI{}{\g},35\SI{}{\W},30\SI{}{\s})$ & $(20\SI{}{\g},20\SI{}{\W},30\SI{}{\s})$\\\hline
    D {\scriptsize(\SI{}{\mm})} & $2.58$ & $4.25$ & $3.69$\\
    \hfill \qquad $\min-\max$ & \checkmark & \checkmark & \checkmark \\ \hline
    W {\scriptsize(\SI{}{\mm})} & $4.88$ & $7.63$ & $7.25$\\
    \hfill \qquad $\min-\max$ & $-11.27\%$ & $-9.17\%$ & \checkmark \\ \hline
    DW {\scriptsize(\SI{}{\mm})} & $0.80$ & $1.47$ & $1.01$\\
    \hfill \qquad $\min-\max$ & \checkmark & $+47\%$ & \checkmark \\ \hline
  \end{tabular}
  \caption{Comparison of the computational lesions against the range of the experimental data.}
  \label{tab_minMaxComparison}
\end{table}

\begin{table}[h!]
  \centering
  \begin{tabular}{l|ccc}
    &\multicolumn{3}{c}{Computational setup ($F$, $P$, $T$)}\\ \hline
    & $(10\SI{}{\g},20\SI{}{\W},30\SI{}{\s})$ &
    $(10\SI{}{\g},35\SI{}{\W},30\SI{}{\s})$ & $(20\SI{}{\g},20\SI{}{\W},30\SI{}{\s})$\\\hline
    D {\scriptsize(\SI{}{\mm})} & $2.58$ & $4.25$ & $3.69$\\
    \hfill \qquad $Q_1 - Q_3$ &  \checkmark & $-1.16\%$ & $+5.42\%$ \\ \hline
    W {\scriptsize(\SI{}{\mm})} & $4.88$ & $7.63$ & $7.25$\\
    \hfill \qquad $Q_1 - Q_3$ & $-22.53\%$ & $-22.93\%$ & $-9.37\%$\\\hline
    DW {\scriptsize(\SI{}{\mm})} & $0.80$ & $1.47$ & $1.01$\\
    \hfill \qquad $Q_1 - Q_3$ & $+14.28\%$ & $+52.22\%$ & $+26.25\%$ \\ \hline
  \end{tabular}
  \caption{Comparison of the computational lesions against the interquartile range $Q_1 - Q_3$ of the experimental data.}
  \label{tab_Q1Q3Comparison}
\end{table}

To dwell deeper into the width underestimation (a feature already observed in \cite{gonzalez2016computational}), we compare the computational lesion dimensions against the interquartile range
$Q_1-Q_3$. The results are shown in Table~\ref{tab_Q1Q3Comparison}. By
comparing Table~\ref{tab_minMaxComparison} and
Table~\ref{tab_Q1Q3Comparison}, it can be seen
that the depth of maximum lesion width is consistently overestimated with respect to the experimental interquartile range, while the lesion width is underestimated. The (\SI{10}{\g}, \SI{35}{\W}, \SI{30}{\s}) setting shows a tiny deviation from the experimental lesion depth, however the experimental results show a very narrow interquartile range of \SI{0.3}{\mm}. 
For the other settings, the depth of the lesion lies in the upper part the interquartile range or is overestimated by $5\%$. Nonetheless, the lesion width is always underestimated while the depth of the maximum width is
overestimated up to \SI{52.22}{\percent}. This shows that the
computational lesion is more spherical than the experimental one.
One possible explanation for this is that our hypothesis of isotropic
heat diffusion in the tissue is not good enough for the type of
prediction we are trying to make.
After all, it is a known fact that cardiac tissue is not isotropic
itself, being made of layers of fibers that run parallel to the surface
of the heart, and rotate clockwise by almost 120 degrees between the endocardium (where the electrode is in contact with the tissue),
and the epicardium \cite{Streeter:1979:GMPG, LeGrice:1995:LSH}.
It is thus worth investigating further the anisotropy of the tissue's thermal  \cite{bhattacharya2003temperature} and electrical conductivities \cite{xie2016effect}. This aspect is currently under study.

\subsection{Comparison of elastic and sharp insertions}
\begin{figure}[h!]
  \centering
  \includegraphics[width=0.3\textwidth]{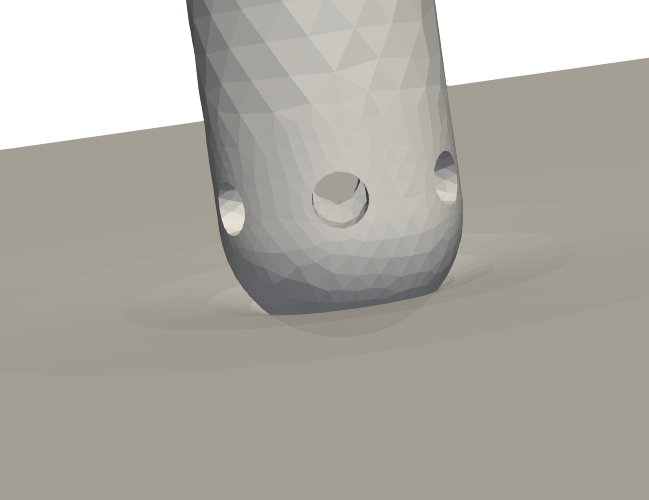}\
  \includegraphics[width=0.3\textwidth]{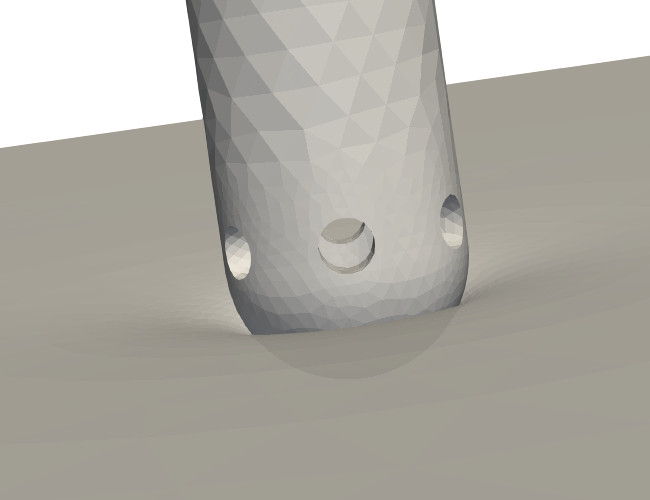}\
  \includegraphics[width=0.3\textwidth]{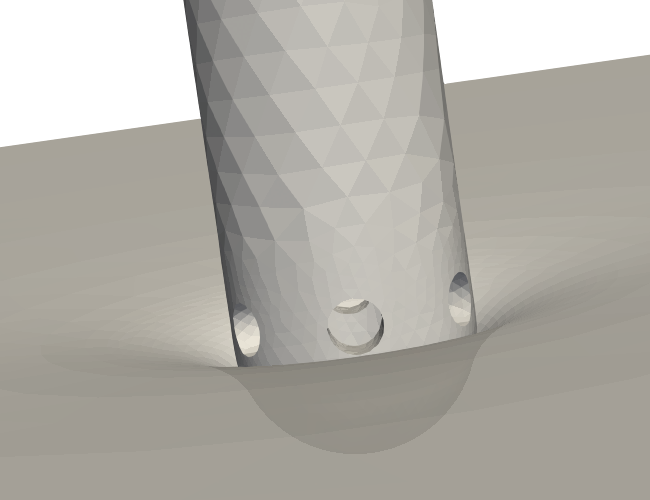}\\
  \vskip0.1cm
  \includegraphics[width=0.3\textwidth]{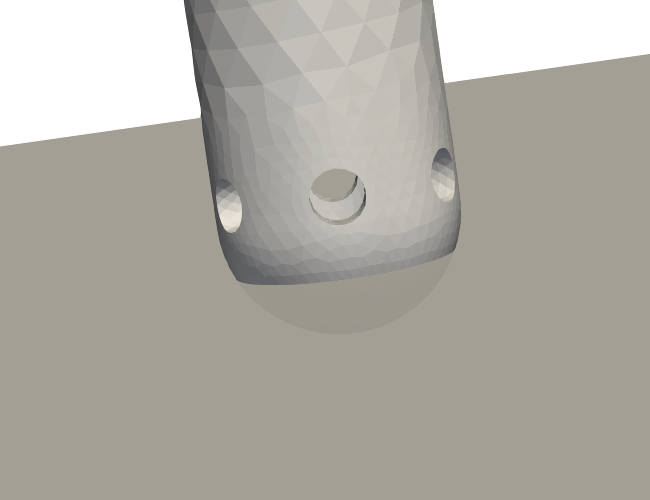}\
  \includegraphics[width=0.3\textwidth]{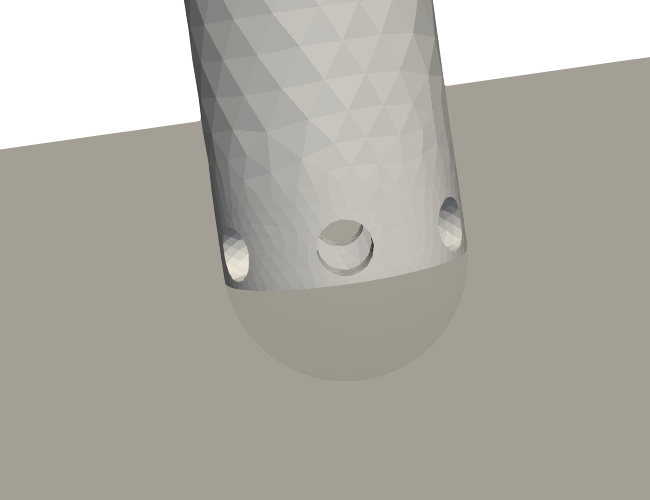}\
  \includegraphics[width=0.3\textwidth]{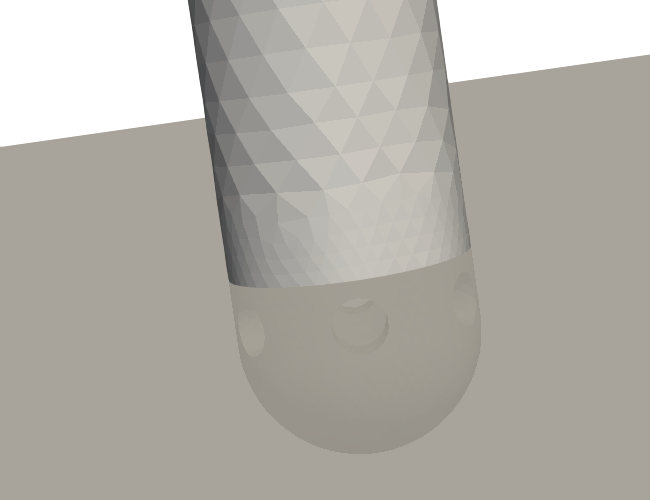}
  \caption{Elastic (top row) and sharp (bottom row) catheter insertion for
  10g (left), 20g (middle), and 40g (right) pressure.}
  \label{fig_comp3}
\end{figure}

One of the strengths (and novelties) of the proposed methodology resides in resorting to an elastic model for the tissue, 
which allows us to treat in a systematic manner both the deformation of the tissue and the 
power dissipated in it. The more accurate description of the tissue-electrode contact surface provided by the elastic model reflects into a more precise evaluation of the power dissipated into the tissue with 
respect to models featuring a sharp insertion of the catheter into an undeformed tissue. The latter is the common configuration encountered in the available literature: besides being pretty unphysical in itself, 
it relies on somewhat arbitrary levels of insertion in the tissue, which are never clearly explained. In addition, a sharp insertion configuration features a larger contact surface between the tissue and the 
electrode. Since the power delivered to the tissue depends on the percentage of electrode surface in contact with the tissue, a sharp insertion is very likely to overestimate the tissue temperature rise and, in consequence, the lesion size. 
To highlight this aspect, we considered two different computational geometries: an elastic deformed tissue and a sharp insertion one. 
In the two geometries, the catheter is placed at the same depth with respect to the undeformed tissue state. 
Figure~\ref{fig_comp3} shows the two different computational geometries for  \SI{10}{\g},  \SI{20}{\g} and  \SI{40}{\g} pressure. 

\begin{figure}[h!]
  \centering
  \includegraphics[width=0.4\textwidth]{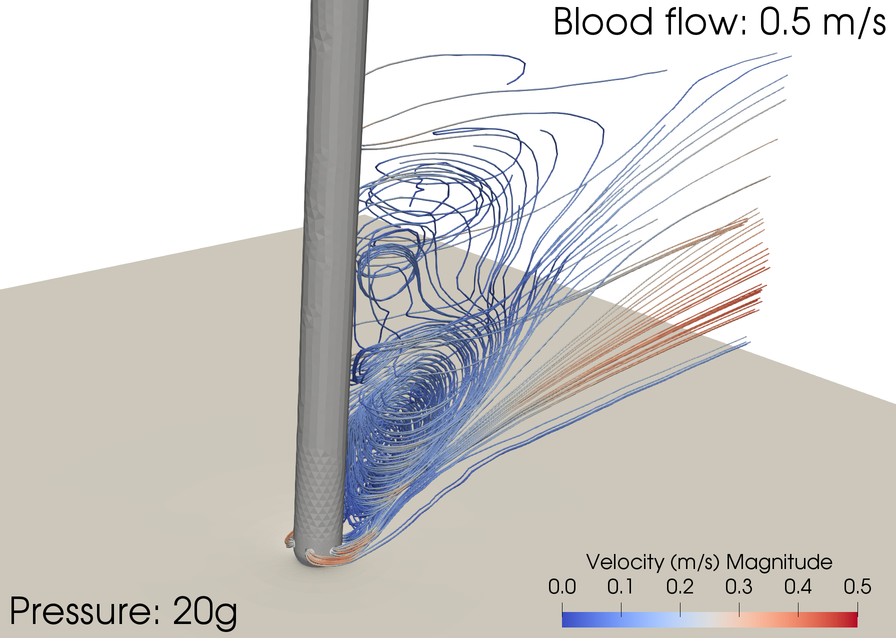}\quad
  \includegraphics[width=0.4\textwidth]{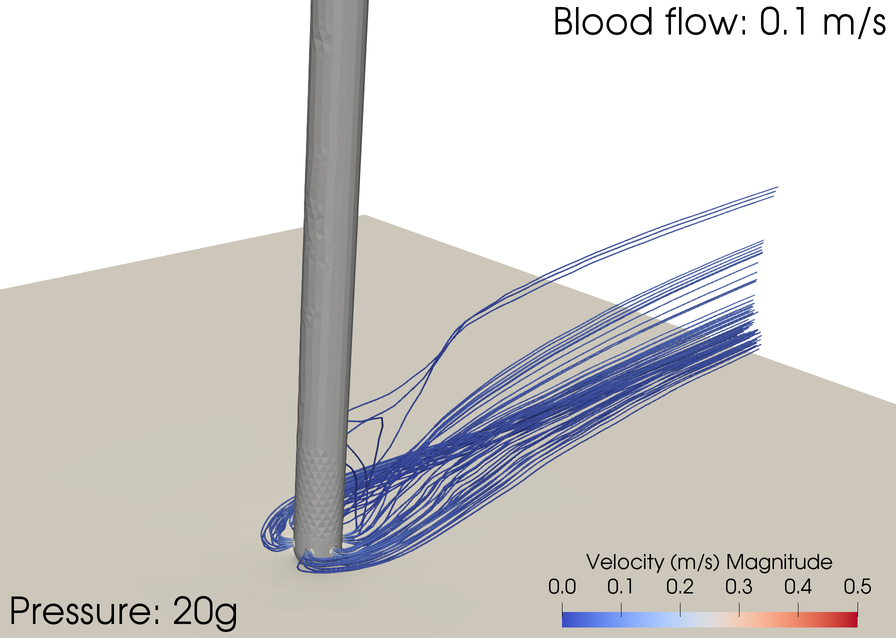}
  \caption{{\bf High-Flow and Low-Flow protocols.} Streamlines of the saline flow from the irrigation holes for the High-Flow (left) and Low-Flow (right) protocols for a catheter pressure of \SI{20}{\g}. }
  \label{fig_saline}
\end{figure}

We simulate an ablation protocol of  \SI{20}{\W} for  \SI{30}{\s}, where we use the tissue parameters of a porcine myocardium and we consider,
inspired from existing experimental studies \cite{guerra2013effects}, two different simulation protocols with respect to the imposed blood flow
velocity: a High-Flow protocol (HF), with $u_b = \SI[per-mode=symbol]{0.5}{\m\per\s}$, and a Low-Flow protocol (LF), with $u_b = \SI[per-mode=symbol]{0.1}{\m\per\s}$.
In Figure~\ref{fig_saline} we show the impact of the two protocols on the streamlines of the saline coolant stemming from the irrigation holes in the elastic \SI{20}{\g} pressure case.

\begin{table}[h!]
  \begin{center}
    \begin{tabular}{l l |c c| c c |c c}
      \multicolumn{2}{c|}{\multirow{2}{*}{\hfill (20W, 30s)}} & \multicolumn{2}{c|}{10g}&\multicolumn{2}{c|}{20g}& \multicolumn{2}{c}{40g} \\
      \cline{3-8}
      &   & elastic& sharp$^*$  & elastic & sharp$^*$ & elastic & sharp$^*$  \\
      \hline
      $\alpha$ &(\%) & 8.46 & 18.87 & 13.29 & 30.73& 19.91& 54.57 \\
      \hline
      $D$ &(\si{\mm}) & 2.58 & 4.16  & 3.69 & 2.61& 4.39&2.07  \\
      $W_x$ &(\si{\mm})&  4.79 & 7.67& 7.00 & 6.61& 9.11& 6.62 \\
      $W_y$ &(\si{\mm})&   4.84& 7.88 &  7.22&  6.71&  9.11& 6.57 \\
      $DW_x$ &(\si{\mm})& 0.80 & 1.52 & 0.91 & 0.11 & 0.93 & -1.12 \\
      $DW_y$ &(\si{\mm})& 0.80 & 1.32 & 1.01 & 0.21 & 0.93 & -1.12 \\
      $V$ &(\si{\mm\cubed}) & 38.28 & 162.6 & 123.55 & 88.33  & 260.15 & 92.44 \\
      $S$ &(\si{\mm\squared})&  0 & 0 &  0& 0 & 0 & 0.63 \\
      $T_{max}$  &(\si{\celsius}) &  61.22 & >100 & 76.95 & >100 & 93.17 & >100 \\
      Pop time &(\si{\s}) & - &  21.6  & - & 5.4 & - & 2.5  \\
      \hline
    \end{tabular}
  \end{center}
  \caption{\label{tab_fluxes_HF} {\bf High-Flow protocol}. Lesion characteristics for sharp and elastic insertion.  $^*$ Simulation
  stopped as a pop: lesion characteristics were assessed at pop time.}
\end{table}

\begin{table}[h!]
  \begin{center}
    \begin{tabular}{l l |c c| c c |c c}
      \multicolumn{2}{c|}{\multirow{2}{*}{\hfill (20W, 30s)}} & \multicolumn{2}{c|}{10g}&\multicolumn{2}{c|}{20g}& \multicolumn{2}{c}{40g} \\
      \cline{3-8}
      &   & elastic& sharp$^*$  & elastic & sharp$^*$ & elastic & sharp$^*$  \\
      \hline
       $\alpha$ &(\%) & 8.46 & 18.87 & 13.29 & 30.73& 19.91& 54.57 \\
     \hline
      $D$ &(\si{\mm}) & 2.46  & 4.17  & 3.62  & 2.63  & 4.38  & 2.09  \\
      $W_x$ &(\si{\mm})& 4.44   & 7.74  & 6.69  & 6.77  & 9.08  & 6.76  \\
      $W_y$ &(\si{\mm})& 4.41  & 7.91  & 6.86  & 6.83  & 9.21  & 6.70  \\
        $DW_x$ &(\si{\mm})& 0.68  & 0.78  & 0.84  & 0.04  & 0.91  & -1.50  \\
        $DW_y$ &(\si{\mm})&  1.03  & 1.33  & 1.04  & -0.26  & 0.91  & -1.20  \\
      $V$ &(\si{\mm\cubed}) &  31.28  & 166.33  & 107.14  & 94.20  & 263.54  & 101.94  \\
      $S$ &(\si{\mm\squared})&  0  & 0  & 0  & 0.07  & 0  & 11.84  \\
      $T_{max}$  &(\si{\celsius}) & 60.43  & >100  & 75.66  & >100  & 92.64  & >100  \\
      Pop time &(\si{\s}) & - & 21.7  & - & 5.3 & - & 2.4  \\
      \hline
    \end{tabular}
  \end{center}
  \caption{\label{tab_fluxes_LF} {\bf Low-Flow protocol}. Lesion characteristics for sharp and elastic insertion. $^*$ Simulation
  stopped as a pop occurred: lesion characteristics were assessed at pop time.}
\end{table}
In Tables~\ref{tab_fluxes_HF} and \ref{tab_fluxes_LF} we report the percentage $\alpha$ of power delivered to the tissue as well as the resulting computational
lesion characteristics for the two insertion profiles in the high blood
flow and low blood flow protocols, respectively.
The tissue power percentage is calculated as described in Section~\ref{sub_pow_calibration}. 
Due to the difference in the contact surface area of the electrode with the tissue, the power delivered to the tissue in the sharp insertion case is more than twice the one of the elastic case. 
Moreover, as this ratio grows with the insertion depth (2.15, 2.33, and 2.41 times for the \SI{10}{\g}, \SI{20}{\g} and \SI{40}{\g} cases, respectively), the tissue temperature rise is expected to be more rapid and overestimated. The lesion depth is assessed as in the previous section. Concerning lesion width and depth of maximal width, we split the measurement in 
two: one along the direction of the blood flow ($W_x$ at depth $DW_x$), the other orthogonal to it ($W_y$ at depth $DW_y$). This choice allows us to assess the difference in lesion anisotropy between the 
two cases. Finally, the lesion volume and surface area are assessed computing the volume enclosed by the
\SI{50}{\celsius} isothermal surface in the tissue, and the area intercepted on the tissue surface by the same isothermal surface.

\begin{figure}
  \centering
  \includegraphics[width=0.4\textwidth]{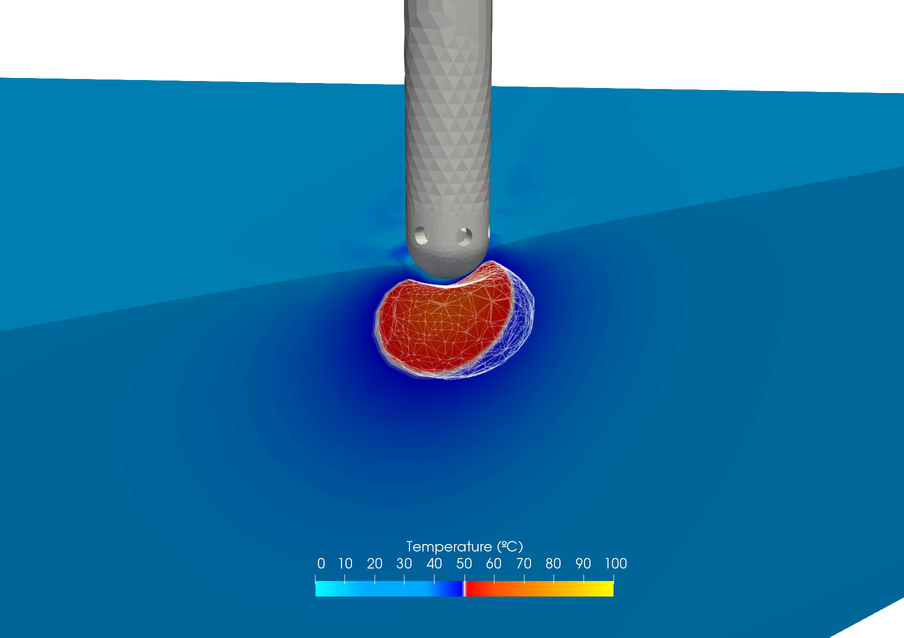}\ \includegraphics[width=0.4\textwidth]{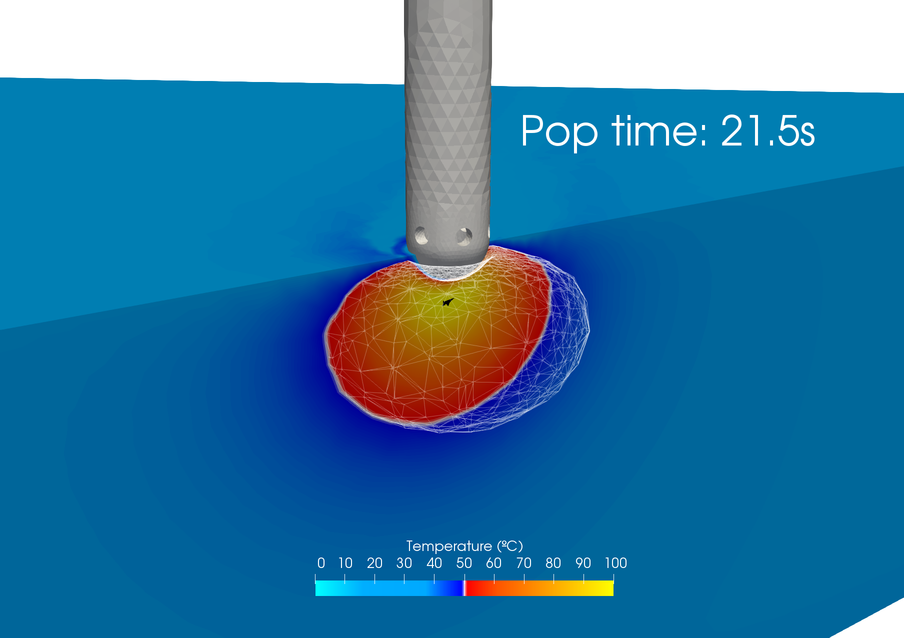}\\ \vskip0.1cm
  \includegraphics[width=0.4\textwidth]{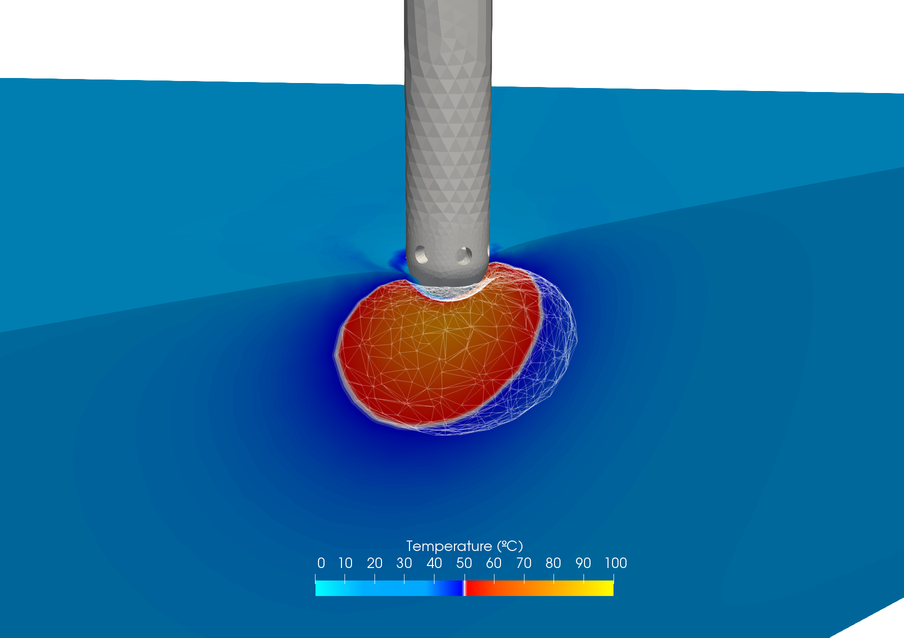}\ \includegraphics[width=0.4\textwidth]{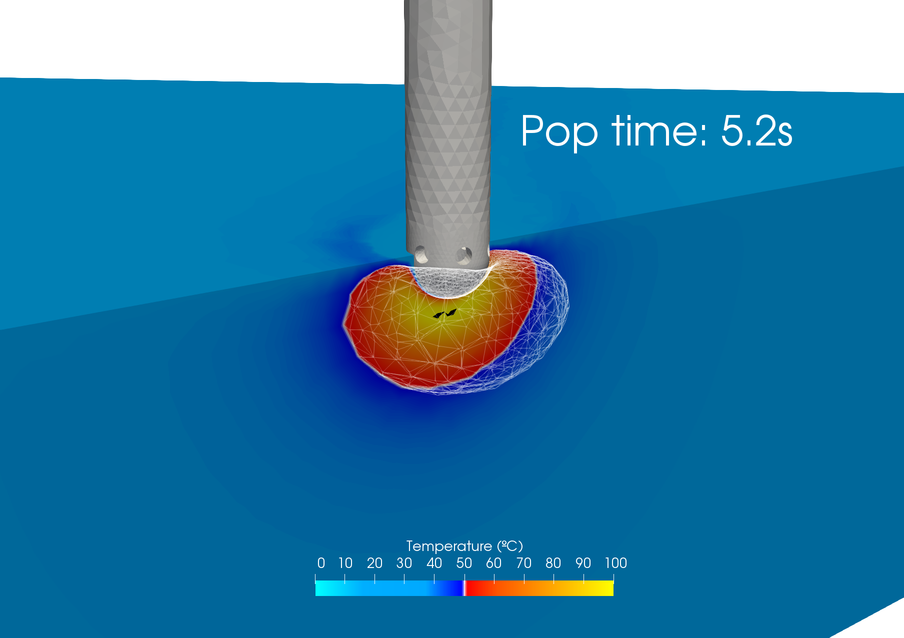}\\  \vskip0.1cm
  \includegraphics[width=0.4\textwidth]{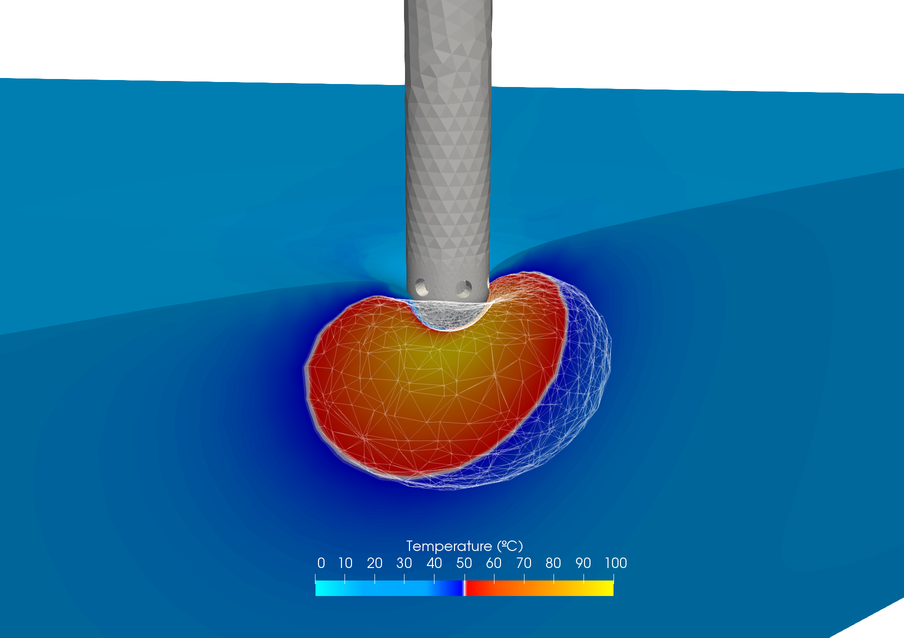}\ \includegraphics[width=0.4\textwidth]{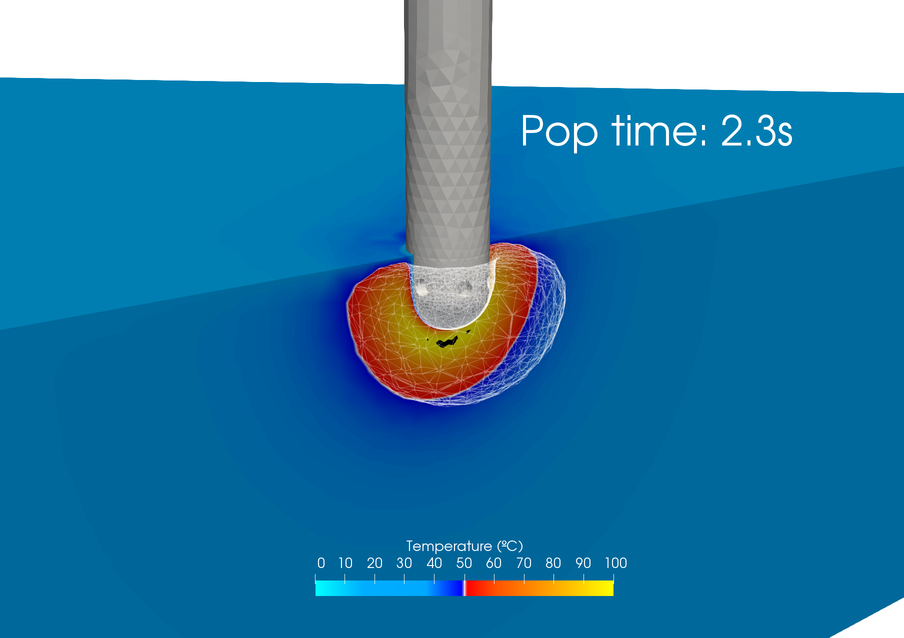}
  \caption{{\bf High-Flow protocol.} Tissue heating with elastic (left column) and sharp (right
  column) catheter insertions for different values of pressure. The 3D lesions are highlighted by the white wireframes. The black spots in the right panels highlight the steam pop location. Top row:
  \SI{10}{\g}. Middle row: \SI{20}{\g}. Bottom row: \SI{40}{\g}. The temperature scale is the
  same for all plots. Simulation protocol: \SI{20}{\W}, \SI{30}{\s}. }
  \label{fig_comp5}
\end{figure}

\begin{figure}
  \centering
  \includegraphics[width=0.4\textwidth]{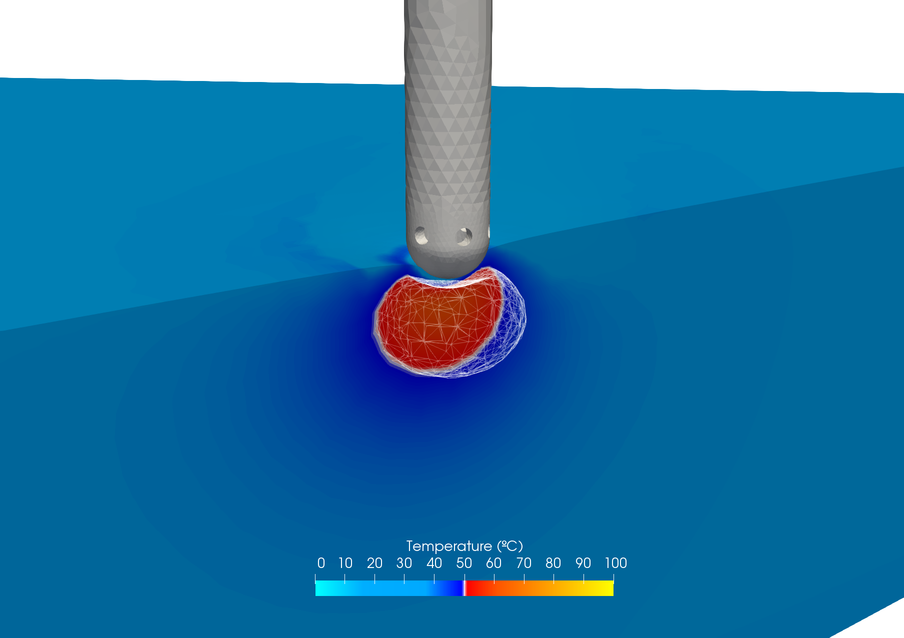}\ \includegraphics[width=0.4\textwidth]{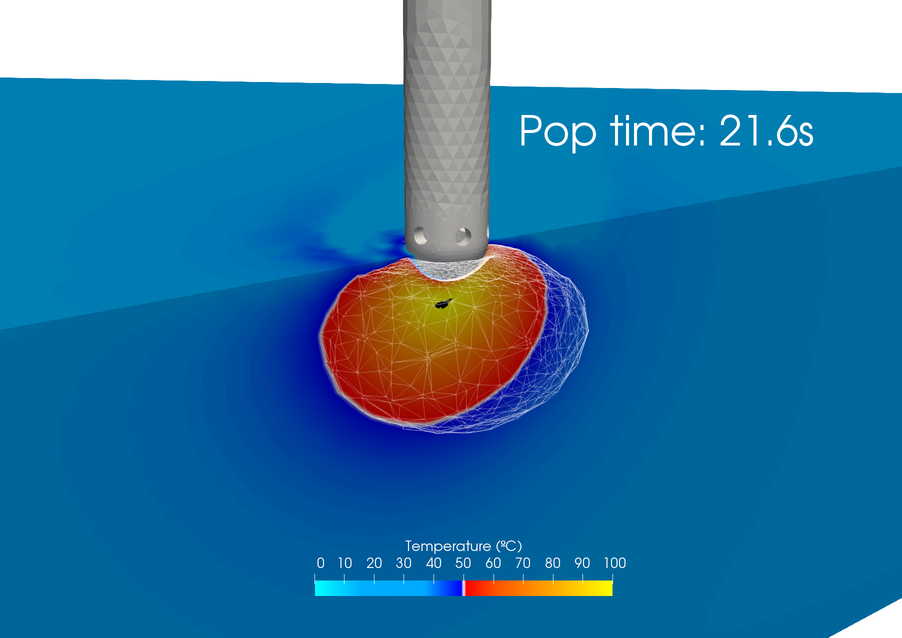}\\ \vskip0.1cm
  \includegraphics[width=0.4\textwidth]{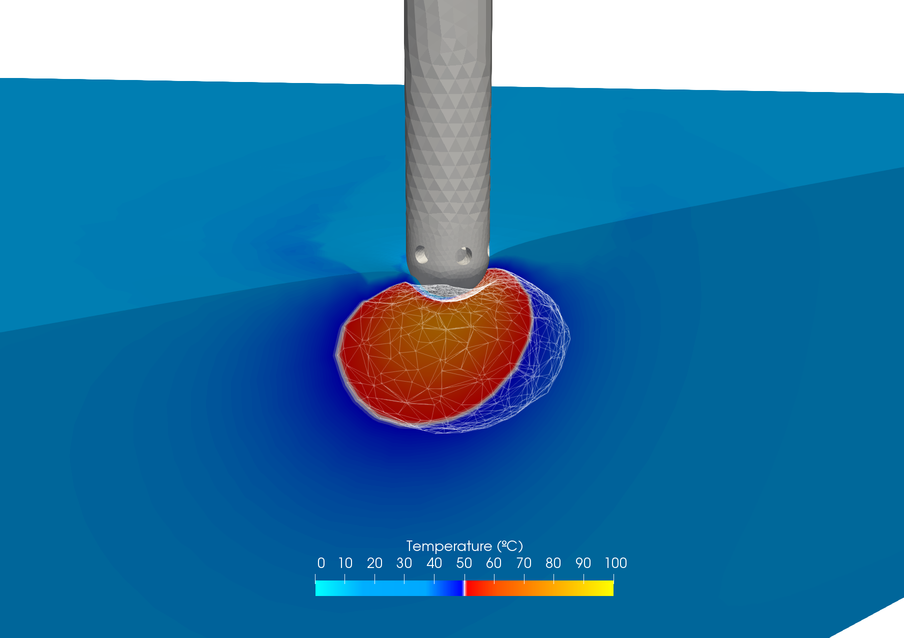}\ \includegraphics[width=0.4\textwidth]{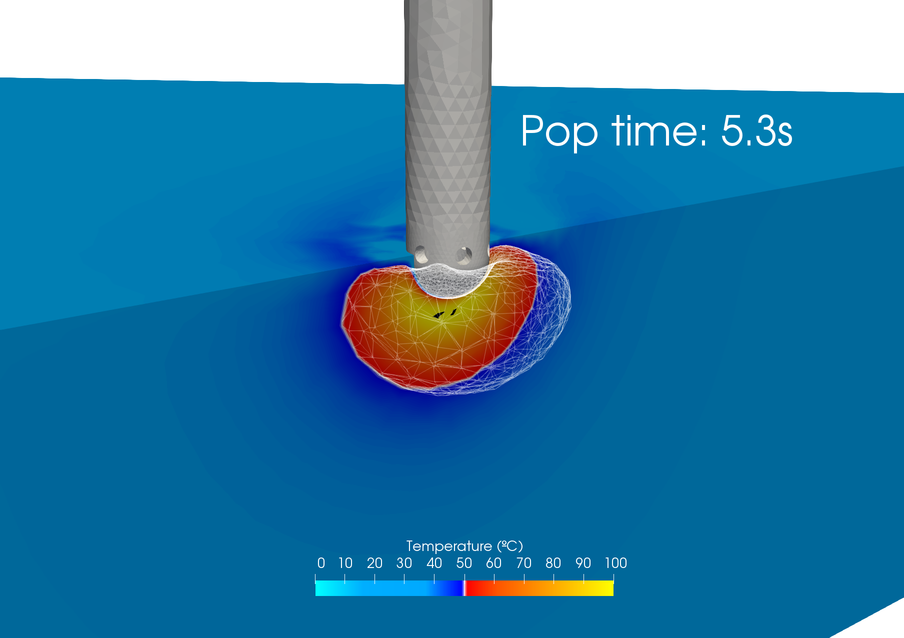}\\  \vskip0.1cm
  \includegraphics[width=0.4\textwidth]{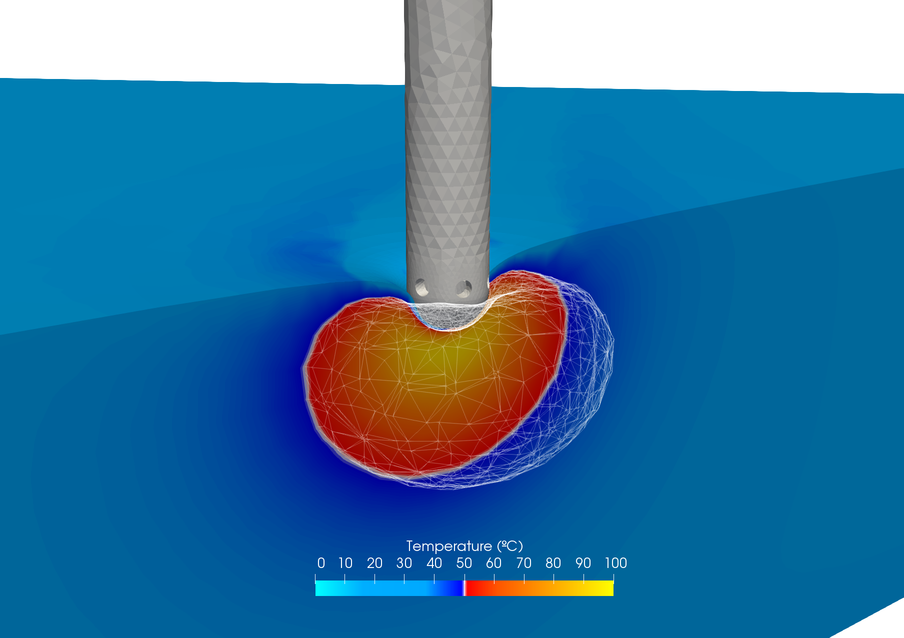}\ \includegraphics[width=0.4\textwidth]{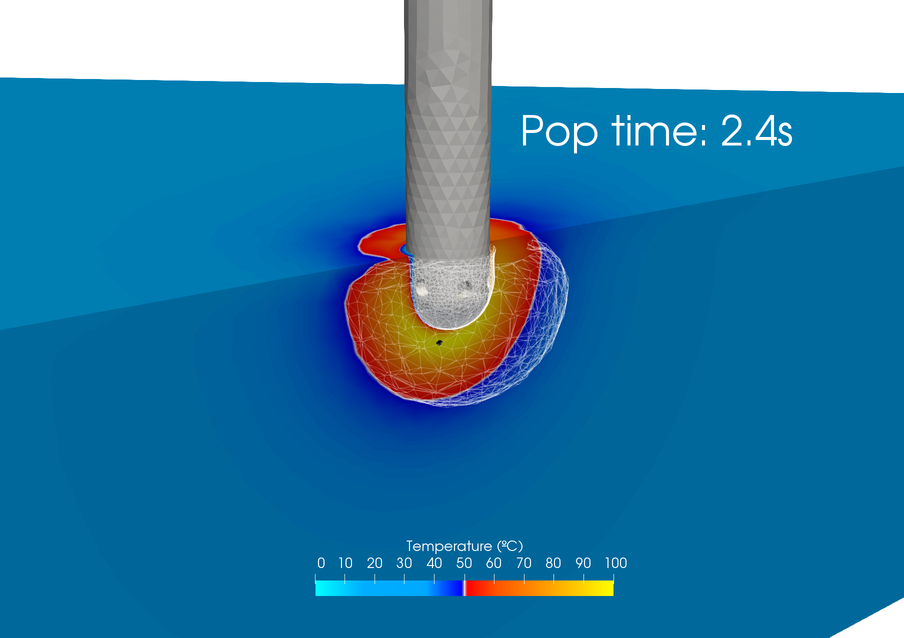}
  \caption{{\bf Low-Flow protocol.} Tissue heating with elastic (left column) and sharp (right
  column) catheter insertions for different values of pressure. The 3D lesions are highlighted by the white wireframes. The black spots in the right panels highlight the steam pop location. Top row:
  \SI{10}{\g}. Middle row: \SI{20}{\g}. Bottom row: \SI{40}{\g}. The temperature scale is the
  same for all plots. Simulation protocol: \SI{20}{\W}, \SI{30}{\s}. }
  \label{fig_comp6}
\end{figure}

From Tables~\ref{tab_fluxes_HF} and \ref{tab_fluxes_LF} we can observe
that the sharp insertion geometry produces consistently larger maximum temperatures and steam pops,
as expected from the bigger percentage of power delivered to the tissue.
In the case of elastic \SI{10}{\g}, the maximum tissue temperature is
approximately \SI{60}{\celsius} for both protocols. On the contrary, a maximum temperature of \SI{100}{\celsius} is reached for the sharp cases after \SI{21.6}{\s}, hinting for the occurrence of steam pops. For both blood flow protocols, the lesion size before popping is much larger than the one developed using an elastic deformation of the tissue. In the \SI{20}{\g} and \SI{40}{\g} cases, the sharp insertion cases have a very rapid temperature increase to \SI{100}{\celsius} in only \SI{5.4}{\s} and \SI{2.5}{\s} respectively, hinting for the occurrence of steam pops. However, a lesion is formed using the elastic tissue and the temperature is below \SI{100}{\celsius} after the completion of \SI{30}{\s}. This comes in good agreement with our \emph{in-vitro} experimental results, where no steam pops were observed in the \SI{20}{\g} and \SI{40}{\g} cases (see Section~\ref{sub_expSetup}). 
In Figures~\ref{fig_comp5} and
\ref{fig_comp6} we compare, for both High and Low-Flow protocols,
the lesions obtained at the end of the simulation (or at the time of popping). The clipping plane for both figures is the diagonal of the computational domain.
The temperature scale is the same for all plots to help visualization. The black spots in the sharp insertion cases highlight the pop location.

\begin{figure}
  \centering
  \includegraphics[width=0.49\textwidth]{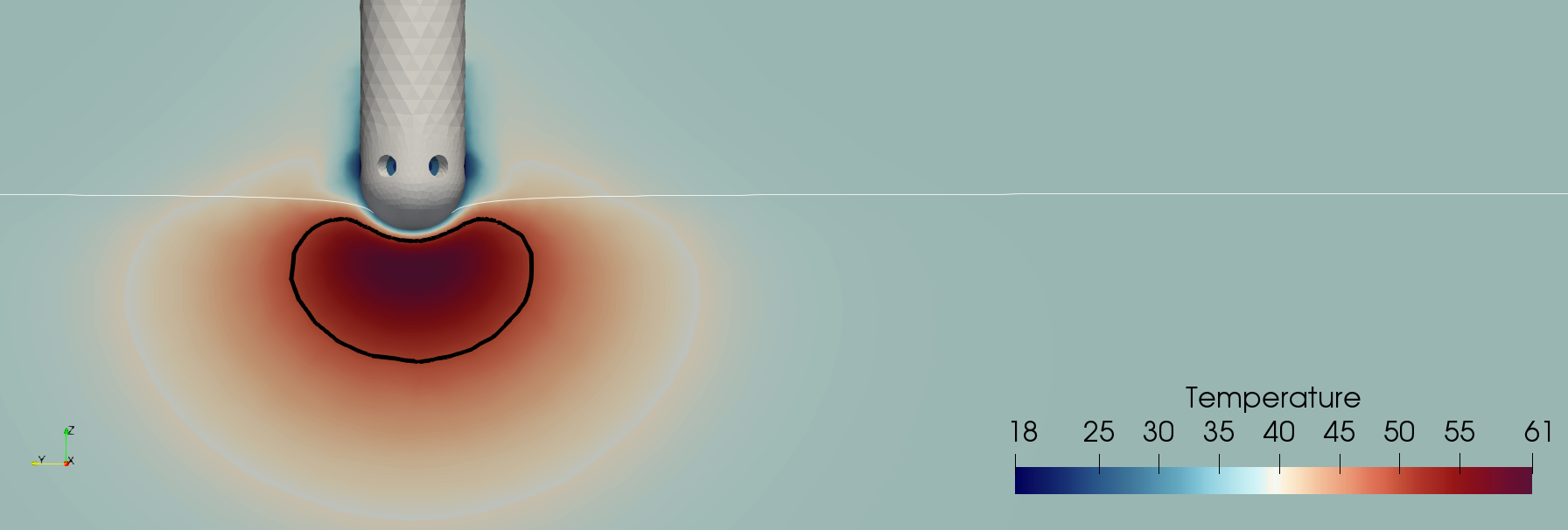}\ \includegraphics[width=0.49\textwidth]{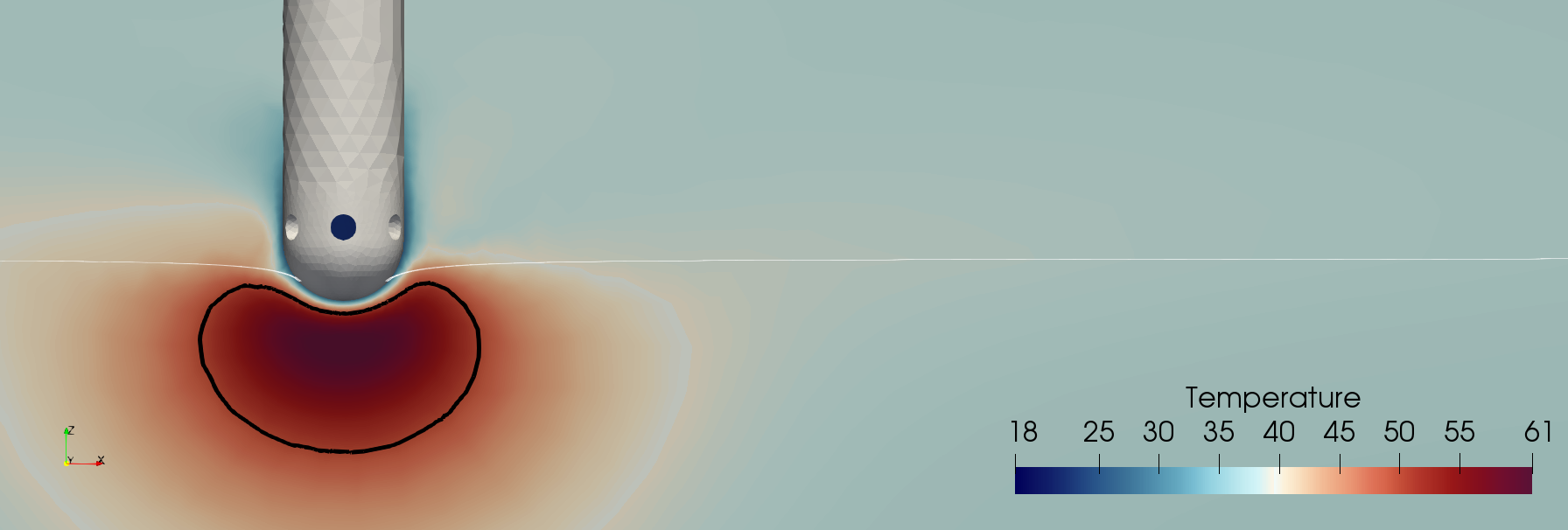}\\ \vskip0.1cm
  \includegraphics[width=0.49\textwidth]{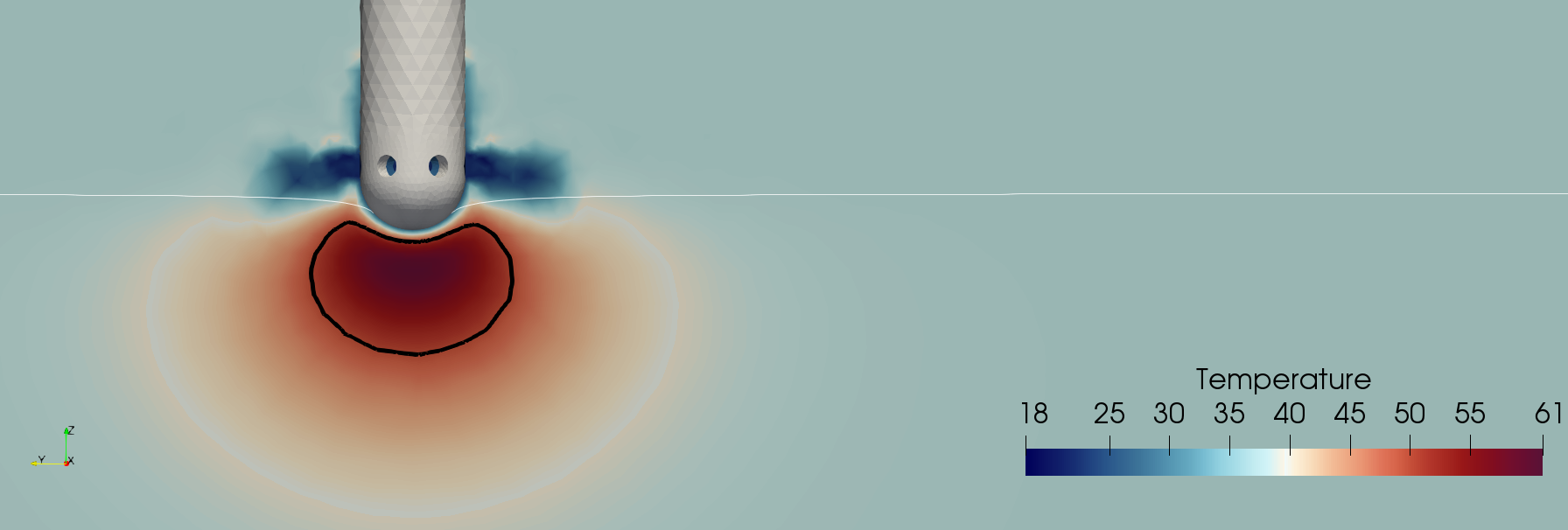}\ \includegraphics[width=0.49\textwidth]{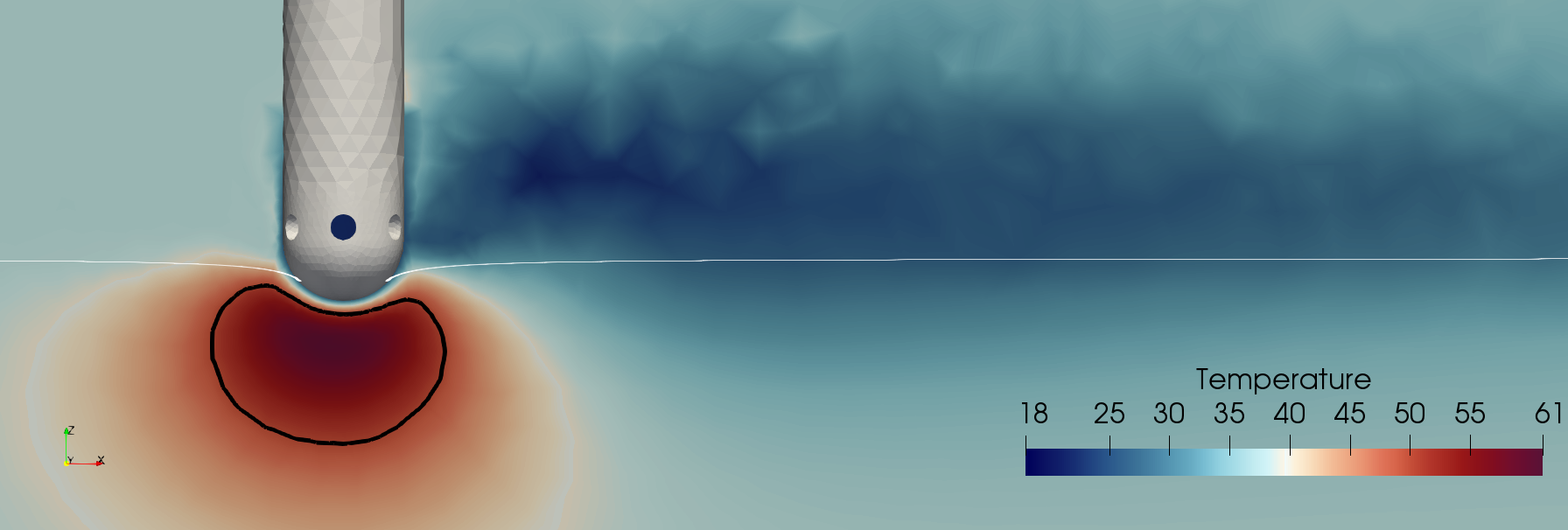}\\  \vskip0.1cm
  \caption{\textbf{High-Flow vs Low-Flow.} The effect of the High-Flow (Top) and Low-Flow (Bottom) protocol in the medial (Left) and lateral (Right) for the elastic catheter insertion case. The computational lesion is displayed in black outline. Simulation protocol: \SI{10}{\g}, \SI{20}{\W}, \SI{30}{\s}.}
  \label{fig_HF_LF}
\end{figure}

Finally, we explore the impact of the high and low blood flow protocols on the elastic cases. It appears that smaller lesions are created in the Low-Flow protocol for \SI{10}{\g} force, while for stronger catheter contact forces such as \SI{40}{\g} the lesion size is similar for both protocols. This effect is a consequence of a larger portion of the catheter being in direct contact with the tissue and the fact that the lesion is formed deeper, thus temperature changes in blood-tissue interface have a negligible effect in the lesion size. For further investigation, we consider the lesion width along the $x$ (along the blood flow)
and $y$ (perpendicular to the blood flow) directions to evaluate the lesion morphological changes.
The High-Flow protocol produces lesions which are more symmetric, with the depth of the maximum width in the medial and lateral direction with respect to the flow being the same. However, differences in the symmetry of the lesion occur for the Low-Flow case.
Figure~\ref{fig_HF_LF} shows the temperature distribution for high and low blood flows for the elastic deformed tissue case and \SI{10}{\g} contact force. Note that the saline cooling effect dominates in the Low-Flow protocol, and leads to a different lesion morphology near the blood-tissue interface. In particular, a small change appears in the width of the lesion in the medial direction with respect to the flow, without affecting the symmetry of the lesion shape. A more evident change in the shape of the lesion appears in the lateral section, where the lesion appears to be tilted opposite to the blood flow due to the strong cooling effect of the irrigated saline. This effect becomes less evident for higher contact force profiles, such as \SI{20}{\g} and \SI{40}{\g}, and the lesion symmetry is restored in both the medial and lateral directions.

\section{Discussion and future work}
\label{sec_discussion}
In this paper we propose a novel framework to model the radiofrequency catheter ablation process, 
which aims at surpassing important limitations of state of the art models.
First, our approach does not only consider the biophysical, but also the mechanical properties
of the cardiac tissue. The model captures the elastic deformation that is produced on the tissue when the catheter exerts a pressure on the cardiac wall. 
An accurate modelling of the catheter geometry allows to account for a more realistic interaction between the coolant liquid springing from the open-irrigated catheter 
and the blood flow. The addition of a board between the cardiac tissue and the dispersive electrode
has a double motivation. On the one hand, neither in {\it in-vitro} experiments nor in clinical settings the 
dispersive electrode is placed in direct contact with the ablated tissue: assuming a direct contact would 
therefore lead to an excessive simplification of the model. On the other hand, the presence of the board allows 
us to use as much information as possible from the RFA system: by properly tuning the board conductivity, 
we can match both the appropriate power delivered to the tissue and the overall system impedance, without 
the need to modify the tissue conductivity. Finally, the amount of power dissipated in the tissue is computed as 
a function of the surface of contact between the electrode and the tissue itself.  
The computational model relies on open-source software, from geometry construction to numerical simulations, 
and postprocessing. In particular, Salome ({\it http://www.salome-platform.org}) was used for mesh generation, 
FEniCS-HPC ({\it http://www.fenics-hpc.org}) for the numerical solution of the PDE system and Paraview 
({\it http://www.paraview.org}) for the postprocessing of the numerical results and the visualizations.

The computational model has been validated in comparison with two sets of 
\emph{in-vitro} experiments, featuring 3 different ablation protocols: (\SI{10}{\g}, \SI{20}{\W}, \SI{30}{\s}),
(\SI{10}{\g}, \SI{35}{\W}, \SI{30}{\s}), (\SI{20}{\g}, \SI{20}{\W} and \SI{30}{\s}).
The numerical lesions have been assessed against the experimental ones through 
some characteristic dimensions such as lesion depth, width and depth of maximum width. 
The simulated lesions are globally in good agreement with the experimental measurements. In particular, the lesion depth always lies within the range of the corresponding experimental results, but the lesion width can be underestimated.
A significant difference is  observed in the \SI{10}{\g} ablation simulation, where the lesion width is 
underestimated by around \SI{10}{\percent}, and in the depth of the maximum width for the \SI{35}{\W} case, which is overestimated by nearly $50\%$.
A more careful comparison with the interquartile $Q_1-Q_3$ of the experimental results 
(see Table~\ref{tab_Q1Q3Comparison}) shows that the depth of maximum lesion width 
is overestimated 
while the lesion width is underestimated. The simulated lesion is thus consistently 
more spherical than the experimental one, prompting for further investigation to clarify this issue. 
This shortcoming can be a consequence of the isotropic heat diffusion coefficient we used 
in this paper: given the lack of isotropy in the cardiac tissue itself, we are currently investigating 
possible anisotropies in the model parameters, in particular the thermal 
\cite{bhattacharya2003temperature} and electrical conductivities \cite{xie2016effect}.

We then compared the effect of using an elastic deformation of the tissue under the catheter pressure 
against profiles featuring a sharp insertion. Accounting for tissue elasticity provides a realistic 
relation between the contact force and the indentation depth in the cardiac tissue. 
As a consequence, the actual power dissipated in the tissue can be computed in a systematic manner, 
differently from other state-of-the-art models featuring sharp-insertion configurations, where it is unclear 
how the same information is deduced.
We compared the two insertion profiles by placing the catheter at the same depth with respect to the undeformed tissue. 
The larger amount of tissue surface in contact with the electrode in the sharp insertion case leads to a systematic overestimation of the temperature rise in the tissue.
To highlight this aspect we simulated an ablation of \SI{20}{\W} for \SI{30}{\s} on a porcine myocardium, with contact forces of \SI{10}{\g}, \SI{20}{\g} and \SI{40}{\g}, and two different blood flow protocols. The sharp insertion cases consistently reach \SI{100}{\celsius} in less than \SI{30}{\s}, which are unphysical taken into account our \emph{in-vitro} experimental results, where no steam pops occurred. Furthermore, in the \SI{40}{\g}
case, the saline irrigation holes are submerged in the cardiac tissue for
a sharp profile, contrary to what has been observed during our \emph{in-vitro} experiments. A visual inspection on the simulated lesions shows that the elastic deformation case produces lesions that 
are closer in shape to the ones obtained in the experiments.

Finally, we compared the effect of the high and low blood flow protocols on the resulting lesions in the elastic case. It appears that the Low-Flow protocol produces smaller lesions for contact forces of \SI{10}{\g} and \SI{20}{\g}, while the lesion size in both protocols are comparable for a contact force of \SI{40}{\g} in both protocols. In addition, we explored the morphological changes of the lesion between the two types of blood flows. For a contact force of \SI{10}{\g}, the Low-Flow protocol affects the lesion orientation, which becomes tilted opposite to the blood flow in the lateral direction due to the cooling of the irrigated saline. This effect becomes less evident as the contact force increases, and for \SI{40}{\g} contact force the lesion symmetric shape is practically restored.

As accurate predictions of the lesions generated by radiofrequency cardiac ablation are fundamental to devise possible new interventional strategies, 
 the use of elastic deformation for the tissue presented in this paper is a significant step towards a more realistic description.  
However, a number of limitations are still present.  First, we considered an isotropic and homogeneous cardiac tissue, even though the
myocardium structure is highly complex and its surface is not smooth
\cite{ho1999anatomy}. In addition, the deformations are assumed to lie within the
elastic bound of the tissue, however it is evident in the experimental
lesions that an elasto-plastic contact occurs. This might be a result of
the temperature change, which alters the mechanical properties of the
tissue, an aspect worth investigating. Finally, in this work we considered is constant Young's modulus of
elasticity, even though a dependency on the heart cycle (systole or diastole) has been observed
\cite{couade2009quantitative}. The effect of the Young's modulus of
elasticity in the proposed RFA model is currently under investigation.


\paragraph{Acknowledgements}
This research was supported by the Basque Government through the BERC 2014-2017 and BERC 2018-2021 program and by Spanish Ministry of Economy and Competitiveness MINECO through BCAM 
Severo Ochoa excellence accreditations SEV-2013-0323 and SEV-2017-0718, and through projects MTM2015-69992-R and MTM2016-76016-R. JJ acknowledges support from project EU H2020 
MSO4SC.  ML acknowledges the "LaCaixa 2016" PhD Grant.  The experimental protocol was supported by Abbott through a non-conditional grant to JG. We acknowledge Esther Jorge-Vizuete for her valuable 
help during the in vitro experiments.


\end{document}